\documentclass[11pt]{amsart}
\usepackage{amsmath,amsfonts,amscd,amssymb}
\theoremstyle{plain}
\newcounter{thmcount}[section]

\newtheorem{theorem}[thmcount]{Theorem}
\newtheorem{conjecture}[thmcount]{Conjecture}
\newtheorem{proposition}[thmcount]{Proposition}
\newtheorem{lemma}[thmcount]{Lemma}
\newtheorem{corollary}[thmcount]{Corollary}
\theoremstyle{definition}
\newtheorem{remark}[thmcount]{Remark}
\newtheorem*{remark*}{Remark}

\newtheorem{definition}[thmcount]{Definition}
\newtheorem*{notation}{Notation}
\newtheorem*{acknowledgements}{Acknowledgements}
\newtheorem{example}[thmcount]{Example}

\def\O{{\mathcal O}}
\def\F{{\mathbb F}}

\def\Q{{\mathbb Q}}
\def\Qb{{\bar{\mathbb Q}}}
\def\Qp{{\mathbb Q}_p}
\def\Qpb{{\bar{\mathbb Q}}_p}
\def\R{{\mathbb R}}
\def\C{{\mathbb C}}
\def\Z{{\mathbb Z}}
\def\Zp{{\mathbb Z}_p}

\def\triv#1{{\mathbf 1}_{#1}}

\def\RC{{\mathcal C}}
\let\rel\Theta
\def\Kbar{\bar{K}}

\def\omegaX{\omega_{\scriptscriptstyle X}}
\def\omegaXbar{\overline{\omega}_{\scriptscriptstyle X}}
\def\omegaY{\omega_{\scriptscriptstyle Y}}

\def\smallcoprod{\raise1.3pt\hbox{$\,\scriptstyle\coprod\,$}}

\def\rkalg#1#2{\rk(#1/#2)}
\def\rksel#1#2#3{\rk_{#3}(#1/#2)}
\def\rkselrep#1#2#3{\rk_{#3}(#1,#2)}
\def\neron#1{\omega_{#1}^o}

\def\newmathop#1{\expandafter\gdef\csname #1\endcsname{\mathop{\rm #1}\nolimits}}
\newmathop{tr}
\newmathop{rk}
\newmathop{ord}
\newmathop{BSD}
\newmathop{Reg}
\newmathop{Gal}
\newmathop{Ind}
\newmathop{Res}
\newmathop{Frob}
\newmathop{GL}
\newmathop{Aut}
\newmathop{coker}
\newmathop{id}
\newmathop{div}
\newmathop{tors}
\newmathop{Hom}
\newmathop{End}
\newmathop{Sel}
\def\ZG#1{\Z_{\scriptscriptstyle\! #1}}
\def\cH{{\mathcal H}}
\def\K{{\mathcal K}}
\def\X{{\mathcal X}}

\def\deton#1#2{\det\bigl(#1\bigm|#2\bigr)}

\def\overarrow#1{\>\>{\buildrel #1 \over \longrightarrow}\>\>}
\def\overlarrow#1{\>\>{\buildrel #1 \over \longleftarrow}\>\>}

\def\bsdc{Birch--Swinnerton-Dyer Conjecture}

\newcommand{\vabove}[2]{\genfrac{}{}{0pt}{}{#1}{#2}}

\def\smallmatrix#1#2#3#4{
  \genfrac{(}{.}{0pt}{1}{#1}{#3}
  \genfrac{.}{)}{0pt}{1}{#2}{#4}
}

\let\lar\longrightarrow
\let\iso\cong
\let\tensor\otimes

\input cyracc.def       % sha
\font\tencyr=wncyr10
\def\sha{\text{\tencyr\cyracc{Sh}}}
\font\eightcyr=wncyr8

\def\smallsha{\text{\eightcyr\cyracc{Sh}}}

\def\beq{$$\begin{array}{llllllllllllllll}}
\def\eeq{\end{array}$$}

\def\ufootnote#1{\insert\footins{\noindent\footnotesize{{\hskip
  1.5em}\llap{${}^{\vphantom a}$}#1}}}

\begin{document}

\title{On the Birch--Swinnerton-Dyer quotients modulo squares}
\author{Tim and Vladimir Dokchitser}
%\date{April 9, 2007}
\address{Robinson College, Cambridge CB3 9AN, United Kingdom}
\ufootnote{2000 {\it Mathematics Subject Classification\/}: Primary 11G40; Secondary 11G05, 11G07, 11G10, 14G25}
%\classification{...}
\email{t.dokchitser@dpmms.cam.ac.uk}
\address{Gonville \& Caius College, Cambridge CB2 1TA, United Kingdom}
\email{v.dokchitser@dpmms.cam.ac.uk}
\dedicatory{to John Coates}
\begin{abstract}
Let $A$ be an abelian variety over a number field $K$.
An identity between the $L$-functions $L(A/K_i,s)$ for extensions
$K_i$ of $K$ induces a conjectural relation between the
Birch--Swinnerton-Dyer quotients.
We prove these relations modulo finiteness of $\smallsha$, and give
an analogous statement for Selmer groups. Based on this, we develop a
method for determining
the parity of various combinations of ranks of $A$
over extensions of~$K$.
As one of the applications, we establish the
parity conjecture for elliptic curves assuming finiteness
of $\smallsha(E/K(E[2]))[6^\infty]$ and some restrictions on the reduction at primes
above 2 and 3: the parity of the Mordell-Weil rank of $E/K$ agrees
with the parity of the analytic rank, as determined by the root number.
We also prove the $p$-parity conjecture for all elliptic curves over $\Q$
and all primes $p$: the parities of the $p^\infty$-Selmer rank
and the analytic rank agree.
\end{abstract}

\llap{.\hskip 10cm} \vskip -0.3cm
\maketitle
\tableofcontents

\section{Introduction}\label{sintro}

The celebrated
conjecture of Birch, Swinnerton-Dyer and Tate asserts
that for every elliptic curve $E$
over a number field $K$,
its Mordell-Weil rank coincides with
the order of vanishing of its $L$-function at $s=1$.
The parity of the latter is determined by the
root number $w(E/K)=\pm 1$, the sign in the expected
functional equation of the $L$-function, leading to

\begin{conjecture}[Parity Conjecture]
\label{parityconj}
The Mordell-Weil rank $\rkalg EK$ is even if and only if the root number
$w(E/K)$ is $+1$.
\end{conjecture}

\noindent
Save for the rank 0 and 1 cases over~$\Q$, virtually nothing is known
about this problem.
At best, one can only lay hands on the
$p^\infty$-Selmer rank $\rksel EKp$ for a prime $p$, that is
the Mordell-Weil rank plus
the number of copies of $\Q_p/\Z_p$ in the Tate-Shafarevich group $\sha(E/K)$.
The Parity Conjecture can also be formulated for Selmer ranks,
as $\sha$ is supposed to be finite by the Shafarevich--Tate conjecture:

\begin{conjecture}[$p$-parity]
\label{pparityconj}
$\rksel EKp$ is even if and only if $w(E/K)=1$.
\end{conjecture}

In view of the conjectures,
the definition of the root number as a product of local terms (local
root numbers) suggests that the parities of $\rkalg EK$ and $\rksel EKp$
should be governed by local data of the elliptic curve.
The purpose of the paper is to develop a theory
that provides such a ``local-to-global'' expression for various
combinations of ranks of $E$ over extensions of~$K$.
The exact description of these ``computable'' combinations
is a curious group-theoretic problem that we have not addressed.
However, there are enough of them to enable us to prove:

\begin{theorem}
\label{ithmparity}
Assuming the Shafarevich--Tate conjecture,
Conjecture \ref{parityconj} holds
over all number fields for elliptic curves with semistable reduction
at primes $v|6$ and not supersingular at $v|2$.
\end{theorem}

\begin{theorem}
\label{ithmpparity}
Conjecture \ref{pparityconj} holds for all $E/\Q$ and all primes $p$.
\end{theorem}

Our starting point is a conjectural formula implied by
Artin formalism for $L$-functions and
the Birch--Swinner\-ton-Dyer conjecture.
As above, fix an elliptic curve
$E/K$ (or a principally polarised abelian variety). Suppose
$L_i, L'_j$ are finite extensions of $K$ such that the
$\Gal(\bar K/K)$-representations
$\bigoplus_i \Ind_{L_i/K} \triv{L_i}$ and
$\bigoplus_j \Ind_{L'_j/K} \triv{L'_j}$ are isomorphic.
Then
$$
  \prod\nolimits_i L(E/L_i,s)=\prod\nolimits_j L(E/L'_j,s)\>,
$$
by Artin formalism.
Ignoring rational squares,
the conjectural expression for the leading terms at $s=1$ leads to
a relation between the regulators and Tamagawa numbers
that we will refer to as the $\square$-Conjecture.
For instance, for semistable elliptic curves it reads
$$
  {\prod\nolimits_i\Reg(E/L_i)c(E/L_i)}\equiv{\prod\nolimits_j\Reg(E/L'_j)c(E/L'_j)}\pmod{\Q^{*2}},
$$
with $c$ the product of local Tamagawa numbers. We will show
that the
$\square$-Conjecture follows from the Shafarevich-Tate conjecture.

\pagebreak

The crucial observation is that the regulators need not cancel by themselves.
It turns out that their quotient can always be expressed through
a combination of Mordell-Weil ranks,
whose parity is therefore determined by local data.
Here is an illustration of how this works in the simplest possible setting,
semistable elliptic curves in $S_3$-extensions:

\begin{example}
Suppose $\Gal(F/K)\iso S_3$,
and let $M,L$ be intermediate extensions of degrees $2$ and $3$ over $K$,
respectively. There is a relation
$$
        (\Ind_{F/K}\triv{F}) \oplus \triv{K}^{\oplus 2}
 \iso   (\Ind_{M/K}\triv{M}) \oplus (\Ind_{L/K}\triv{L})^{\oplus 2}\>.
$$
(i) For semistable $E/K$, the $\square$-Conjecture implies that
$$
  \frac{\Reg(E/F)\Reg(E/K)^2}{\Reg(E/M)\Reg(E/L)^2} \equiv
  \frac{c(E/F)c(E/K)^2}{c(E/M)c(E/L)^2} \pmod{\Q^{*2}}\>.
$$
(ii) The quotient of regulators is related to Mordell-Weil ranks
(Ex. \ref{exs3}):
$$
 3^{\rkalg EK+\rkalg EM+\rkalg EL}
    \equiv
 \frac{\Reg(E/F)\Reg(E/K)^2}{\Reg(E/M)\Reg(E/L)^2}
  \pmod{\Q^{*2}}\>.
$$
Thus, assuming finiteness of $\sha$, we obtain an expression
for the sum of the three ranks $\rkalg EK+\rkalg EM+\rkalg EL$
in terms of local data.

\noindent
(ii$'$) In fact, by a somewhat more sophisticated
technique, we can prove an analogous (unconjectural) statement about
$3^\infty$-Selmer ranks (Thm. \ref{thmselft}):
$$
 \rksel EK3+\rksel EM3 +\rksel EL3
  \equiv \ord_3 \frac{c(E/F)c(E/K)^2}{c(E/M)c(E/L)^2} \pmod2\>.
$$
(iii) Finally, a purely local computation allows us to
relate the Tamagawa numbers to root numbers (Prop.~\ref{thmcv}):
$$
  w(E/K)w(E/M)w(E/L)=1 \>\>
    \Longleftrightarrow
  \>\>
  \ord_3\frac{c(E/F)c(E/K)^2}{c(E/M)c(E/L)^2}
  \equiv0\pmod2\>,
$$
and we obtain a special case of the parity conjecture for
$S_3$-extensions.
\end{example}

The layout of the paper is as follows:

In \S\S\ref{ssBSD}--\ref{ssSQUARE} we formulate the $\square$-Conjecture
and prove it assuming finiteness of $\sha$
(Conj. \ref{C2}, Cor. \ref{C2true}).
This relies on invariance of the BSD-quotient under Weil restriction
of scalars and under isogenies.
Next, we relate the quotient of regulators from the conjecture
to the parity of Mordell-Weil ranks in \S\ref{ssmoreregs}
(Thm. \ref{regconst}, Cor. \ref{correg}),
and give examples in \S\ref{ssregex}.

Thus, we have now complete versions of steps (i) and (ii) of the above
example (principally polarised abelian varieties and arbitrary
field extensions). We do not attempt to deal with (iii)
in such generality,
but confine ourselves to elliptic curves and
extensions with Galois group
$\smallmatrix1*0* \subset \GL_2(\F_p)$.
After reviewing the classification of root numbers in \S\ref{ssroot}, we
relate the Tamagawa numbers to root numbers for
such extensions in~\S\ref{sstamagawa} (Prop.~\ref{thmcv}).
Combined with the results of~\cite{TV-P} on the parity conjecture for
elliptic curves with a 2-isogeny, this proves
Theorem~\ref{ithmparity}.

So far, we related parities of Mordell-Weil ranks to Tamagawa numbers
assuming that $\sha$ is finite. In \S\ref{sselmer}
we address the problem of getting an unconditional statement
about Selmer ranks (as in (ii$'$)).
We prove an analogue of the $\square$-Conjecture
(Thm.~\ref{thmsel}, Cor.~\ref{corsel})
by tweaking Tate--Milne's proof of the isogeny invariance of
the \bsdc.
The quotient of regulators becomes replaced by a
quantity $Q$ measuring the effect of an isogeny on Selmer groups.
In \S\ref{ssselmain} we turn $Q$ into Selmer ranks
in fair generality (Thm.~\ref{selmain}, Cor.~\ref{selav}),
and we illustrate it for
$S_n$-extensions (Ex.~\ref{selsym}),
$\smallmatrix1*0*$-extensions (\S\ref{ssselft}),
and dihedral extensions (\S\ref{ssseldih}).
In \S\ref{ssselft} we give an application
to ranks of elliptic curves in false Tate curve towers.
We end in \S\ref{sspparity}
by proving Theorem \ref{ithmpparity}; for odd $p$ it is a consequence of
our results for dihedral extensions and the existence of quadratic
and anticyclotomic twists for which the Birch--Swinnerton-Dyer rank formula
is known to hold.%
\footnote{Since writing of this paper, we have extended (ii') and (iii)
to arbitrary $\Gal(F/K)$ in \cite{Selfduality,Tamroot};
for (ii') the theory is now as clean as for (ii), e.g. computations with
isogenies in \S\S\ref{ssselft}--\ref{ssseldih} are replaced by elementary
representation theory, as in \S\ref{ssregex}.}

\smallskip

Finally, let us mention how the applications of our theory connect to
earlier work.
To our best knowledge, over number fields Theorem \ref{ithmparity}
is the first general result of this kind,
except for the work \cite{CFKS,TV-P} on curves with a $p$-isogeny.
In contrast, the $p$-parity conjecture over $\Q$ was known in almost all
cases, thanks to Birch, Stephens, Greenberg and Guo~\cite{BS,Gre,Guo}~($E$ CM),
Kramer, Monsky \cite{Kra,Mon} ($p=2$), Nekov\'a\v r~\cite{NekS} ($p$ potentially ordinary
or potentially multiplicative) and Kim~\cite{Kim} ($p$ supersingular).
The results for Selmer groups in dihedral and false Tate curve
extensions are similar to those recently obtained by Mazur--Rubin \cite{MR}
and Coates--Fukaya--Kato--Sujatha \cite{CFKS,CS}, respectively.

\begin{notation}
Throughout the paper $K$ always denotes a number field.
For a place $v$ of $K$ we write $|\cdot|_v$ for the normalised absolute value
at $v$.
If $L/K$ is a finite extension, we denote by $\Ind_{L/K}\triv{L}$ the
induction of the trivial (complex)
representation of $\Gal(\Kbar/L)$ to $\Gal(\Kbar/K)$.
This is the permutation representation corresponding to the
set of $K$-embeddings~of~$L$~into~$\Kbar$.

For an elliptic curve $E/K$ we use the following notation:
\smallskip

\begin{tabular}{ll}
$\rkalg EK$           & Mordell-Weil rank of $E/K$. \cr
$\rksel EKp$          & $p^{\infty}$-Selmer rank of $E/K$, i.e. \cr
                      & $\rkalg EK+$ number of copies of $\Q_p/\Z_p$ in $\sha(E/K)$. \cr
$w(E/K_v)$            & local root number of $E$ at a place $v$ of $K$.\cr
$w(E/K)$              & global root number, $=\prod_{v} w(E/K_v)$.\cr
$\Reg(E/K)$           & regulator of $E/K$, i.e. $|\det|$ of the canonical\cr
                      & height pairing on a basis of $E(K)/E(K)_{\tors}$.\cr
$c_v$                 & local Tamagawa number at a finite place $v$.\cr
$c(E/K)$              & product of the local Tamagawa numbers, $=\prod_{v\nmid\infty} c_v$.\cr
$W_{F/K}(E)$          & the Weil restriction of scalars of $E/F$ to $K$.\cr
\end{tabular}
\smallskip

\noindent Finally, we will need a slight modification of $c(E/K)$.
Fix an invariant differential $\omega$ on $E$.
Let $\neron{v}$ be N\'eron differentials at finite places $v$ of~$K$,
and set
$$
  C(E/K) = \prod_{v\nmid\infty} c_v \Bigl|\frac{\omega}{\neron{v}}\Bigr|_{_v}\>.
$$
Note that $C(E/K)$ depends on the choice of $\omega$,
although we have omitted this from the notation.
When writing $C(E/L_i)$ for various extensions $L_i/K$, we always implicitly
use the same $K$-rational differential.

We use similar notation for abelian varieties (the analogue of an invariant
differential being a non-zero global exterior form). We write $A^t$ for the
dual abelian variety.

\begin{acknowledgements}
Parts of this research have been carried out while the first author (T.)
was a University Research Fellowship holder of the Royal Society and while
the second author (V.) stayed at the Max Planck Institute
for Mathematics (Bonn). We are grateful to both institutions for
their support.
We would also like to thank
Henri Darmon,
Claus Diem, % Weil restrictions
Matteo Longo and
Robert Prince % doom music
for helpful discussions,
and the referee for the comments on the paper.
\end{acknowledgements}

\end{notation}

%%%%%%%%%%%%%%%%%%%%%%%%%%%%%%%%%%%%%%%%%%%%%%%
\section{$\square$-Conjecture and regulator quotients}\label{ssquareness}

\subsection{Artin formalism and BSD-quotients}
\label{ssBSD}

Let $K$ be a number field and
let $A/K$ be an abelian variety with a fixed
global exterior form $\omega$. Recall the statement
of the Birch--Swinnerton-Dyer conjecture:

\begin{conjecture}[Birch--Swinnerton-Dyer, Tate \cite{TatC}]
\label{bsd}
\
\begin{enumerate}
\item The $L$-function $L(A/K,s)$ has an analytic continuation to $s=1$,~and
$$
  \ord_{s=1} L(A/K,s) = \rkalg AK \>.
$$
\item The Tate-Shafarevich group $\sha(A/K)$ is finite,
and the leading coefficient of $L(A/K,s)$ at $s=1$ is
\begingroup
$$
\font\txtfnt=cmr9\textfont0=\txtfnt
\font\mthfnt=cmmi9\textfont1=\mthfnt
\font\scrfnt=cmr6\scriptfont0=\scrfnt
\font\sscfnt=cmmi6\scriptfont1=\sscfnt
  \BSD(A/K)=
  \frac{|\smallsha(A/K)|\Reg(A/K)\,C(A/K)}{|A(K)_{\tors}||A^t(K)_{\tors}||\Delta_K|^{\dim A/2}}
  \prod_{\vabove{v|\infty}{\text{real}}} \int\limits_{A(K_v)}\!\!\!\! |\omega|
  \prod_{\vabove{v|\infty}{\text{cplx}}} 2\!\!\!\int\limits_{A(K_v)}\!\!\!\! \omega\wedge \bar\omega.
$$
\endgroup
\end{enumerate}
\end{conjecture}

\begin{notation}
We call $\BSD(A/K)$
the Birch--Swinnerton-Dyer quotient for $A/K$. We also write $\BSD_p(A/K)$
for the same expression with $\sha$ replaced by its $p$-primary component
$\sha[p^\infty]$.
(They are independent of the choice of $\omega$ by the product formula.)
\end{notation}

Now let $L_i\supset K$ and $L'_j\supset K$ be number fields such that
$$
  \bigoplus\nolimits_i \Ind_{L_i/K} \triv{L_i} \iso
  \bigoplus\nolimits_j \Ind_{L'_j/K} \triv{L'_j}
$$
as complex representations of $\Gal(\Kbar/K)$.
In other words, the $\Gal(\Kbar/K)$-sets
$\coprod_i\Hom_K(L_i,\Kbar)$ and $\coprod_j\Hom_K(L'_j,\Kbar)$
give rise to isomorphic permutation representations.
By Artin formalism for $L$-functions,
$$
  \prod\nolimits_i L(A/L_i,s) = \prod\nolimits_j L(A/L'_j,s) \>,
$$
so the following is a consequence of Conjecture \ref{bsd}.

\begin{conjecture}
\label{C1}
With $A/K$ and $L_i, L'_j$ as above,
\begin{itemize}
\item[(a)]
$\sum_i \rkalg A{L_i} = \sum_j \rkalg A{L'_j}$,
\item[(b)]
$\sha(A/L_i), \sha(A/L'_j)$ are finite, and
$
  \prod_i \BSD(A/L_i) = \prod_j \BSD(A/L'_j) \>.
$
\end{itemize}
\end{conjecture}

This is in effect a compatibility statement of the
Birch--Swinnerton-Dyer conjecture with Artin formalism.
Part (a) is easily seen to be true: let $F/K$ be a finite Galois
extension containing $L_i$ and $L'_j$, and let $V=A(F)\tensor_\Z\C$. Then
$$
  \rkalg A{L_i} = \dim V^{\Gal(F/L_i)}
   = \langle \triv{L_i}, \Res_{F/L_i} V \rangle
   = \langle \Ind_{F/L_i}\triv{L_i}, V \rangle
$$
by Frobenius reciprocity, and similarly for $L'_j$; now take the sum over
$i$ and~$j$.

We now show that (b) is implied by finiteness of $\sha$.
As C. S. Dalawat, K. Rubin and M. Shuter pointed out to us, this is
essentially the same as H. Yu's Theorem 5 in \cite{Yu}.

\begin{theorem}
\label{main}
Let $A/K$ be an abelian variety, and let $L_i, L'_j$ be
finite extensions of $K$ satisfying
$\oplus_i \Ind_{L_i/K} \triv{L_i} \iso\oplus_j \Ind_{L'_j/K} \triv{L'_j}$.
Suppose that $\sha(A/L_i), \sha(A/L'_j)$ are finite.
Then Conjecture \ref{C1}b holds.

Furthermore, if we weaken the assumption to
$\sha(A/L_i)[p^\infty], \sha(A/L'_j)[p^\infty]$ being finite for some prime $p$,
then the $p$-part of Conjecture \ref{C1}b holds, i.e.
$$
  \prod_i \BSD_p(A/L_i) \Bigm/ \prod_j \BSD_p(A/L'_j)
$$
is a rational number with trivial $p$-valuation.
\end{theorem}

\begin{proof}
For $F/K$ finite, write $W_{F/K}(A)$ for the Weil restriction of scalars
of $A/F$ to $K$. This is an abelian variety over $K$ of dimension
$[F\!:\!K]\dim A$,~and $\BSD_p(W_{F/K}(A))\!=\!\BSD_p(A/F)$ provided that
$\sha(A/F)[p^\infty]$ is finite~(\cite{MilO}~\S1).

Consider
$$
  X=\prod_i W_{L_i/K}(A), \qquad Y=\prod_j W_{L'_j/K}(A)\>.
$$
Then $\BSD_p(X)=\prod_i\BSD_p(A/L_i)$ and $\BSD_p(Y)=\prod_j\BSD_p(A/L'_j)$.
By the invariance of the Birch--Swinnerton-Dyer quotient under
isogenies (\cite{CasVIII}, \cite{TatC} and \cite{MilA} Thm. 7.3, Remark 7.4),
it suffices to show that $X$ and $Y$ are isogenous.

As the representations
$\oplus_i \Ind_{L_i/K} \triv{L_i}$ and $\oplus_j \Ind_{L'_j/K} \triv{L'_j}$
are realisable over $\Q$ and are isomorphic over $\C$, they are isomorphic
over $\Q$ (see e.g. \cite{SerLi}, Ch. 12, Prop. 33 and remark following it).
So the corresponding integral permutation modules are isogenous, in the
sense that there is an inclusion of one as a finite index submodule of the
other. This induces an isogeny $X\to Y$ (see \cite{MilO}, \S2, Prop. 6a).
\end{proof}

\subsection{$\square$-Conjecture}
\label{ssSQUARE}
Although
Conjecture \ref{C1}b has the advantage that it does not involve $L$-functions,
it still relies on finiteness of $\sha$.
Also, even when $\sha$
is finite it is hard to determine, which makes the statement difficult
to work with. However, if $A$ is principally polarised and $\sha$ is finite,
then the order of $\sha$ is either a square or twice a square by the non-degeneracy of
the Cassels--Tate pairing \cite{TatD}. (If $A$ is an elliptic curve or has
a principal polarisation arising from a $K$-rational divisor, then
the order of $\sha$ is a square, see \cite{CasIV,TatD,PS}.) So we can
eliminate $\sha$ from the statement by working modulo squares, which
also removes the contribution from the torsion. Moreover,
in this combination of BSD-quotients, the discriminants
of fields cancel by the conductor-discriminant formula, as do
the real and complex periods, provided that one chooses
the same $\omega$ over $K$ for each term.
Thus Conjecture \ref{C1}b implies the following
(see Remark \ref{genC2} for an extension to abelian varieties):

\begin{conjecture}[$\square$-Conjecture]
\label{C2}
Let $E/K$ be an elliptic curve, and fix
an invariant differential $\omega$ on $E$. Let $L_i, L'_j$ be
finite extensions of $K$ satisfying
$\oplus_i \Ind_{L_i/K} \triv{L_i} \iso\oplus_j \Ind_{L'_j/K} \triv{L'_j}$.
Then
$$
  \prod_i \Reg(E/L_i)\,C(E/L_i) \equiv \prod_j \Reg(E/L'_j)\,C(E/L'_j) \pmod{\Q^{*2}}.
$$
\end{conjecture}

\begin{corollary}[of Theorem \ref{main}]
\label{C2true}
The $p$-part of Conjecture \ref{C2} holds, provided
that $\sha(E/L_i)[p^\infty]$ and $\sha(E/L'_j)[p^\infty]$ are finite.
In other words,
$$
  \prod_i \Reg(E/L_i)\,C(E/L_i) \Bigm/ \prod_j \Reg(E/L'_j)\,C(E/L'_j)
$$
is a rational number with even $p$-valuation.
\end{corollary}

We are going to explore the (surprisingly non-trivial) consequences of this
for parities of Mordell-Weil ranks. Here is a simple example:

\begin{example}\label{ex2}
Take the modular curve $E=X_1(11)$ over the fields $\Q$,
$\Q(\mu_3)$,  $L=\Q(\sqrt[3]{m})$ and $F=\Q(\mu_3,\sqrt[3]{m})$
for $m>1$ cube free.
We have an equality of $\Gal(\Qb/\Q)$-representations,
$$
  (\Ind_{F/\Q}\triv{F}) \oplus (\triv{\Q})^{\oplus 2}
  \iso (\Ind_{L/\Q}\triv{L})^{\oplus 2} \oplus (\Ind_{\Q(\mu_3)/\Q}\triv{\Q(\mu_3)})\>.
$$
The Mordell-Weil rank of $E$ is 0 over $\Q(\mu_3)$, so
$\Reg(E/\Q)$ and $\Reg(E/\Q(\mu_3))$ are both $1$. The
$\square$-Conjecture implies that
$$
 \frac{\Reg(E/F)}{\Reg(E/L)^2}
 \equiv
 \frac{c(E/\Q(\mu_3))\,c(E/L)^2}{c(E/F)\,c(E/\Q)^2}
 \pmod{\Q^{*2}}\>.
$$
For $v|11$, the local Tamagawa number $c_v$ is the
valuation of the minimal discriminant ($=-11$) at~$v$, because $E$ has split
multiplicative reduction at~$v$. So, by a simple computation,
the above quotient of Tamagawa numbers is $1$ when $11\nmid m$
and $3$ when $11|m$.
On the other hand, let $P_1,..,P_n$ be a basis for $E(L)\tensor\Q$,
and $H$ the height matrix $\langle P_i,P_j \rangle_{L}$, so
that $\Reg(E/L)$ is $|\det(H)|$ up to a (rational) square.
If $g\in\Gal(F/\Q)$ is an element of order $3$, then
$P_1,..,P_n,P_1^g,..,P_n^g$ is a basis for $E(F)\tensor\Q$.
One readily verifies that the height matrix over $F$ is
$\smallmatrix{2H}{-H}{-H}{2H}$, so the regulator $\Reg(E/F)$ is
$3^n|\det(H)|^2$ up to a square.
Hence the $\square$-Conjecture implies that $\rkalg E{\Q(\sqrt[3]{m})}$
is odd if and only if $11|m$.
(See Example \ref{exs3} and Corollary \ref{rkmwft} for a generalisation.)
\end{example}

We end with a few observations:

\begin{remark}\label{genC2}
There are obvious analogues of the
$\square$-Conjecture and Corollary \ref{C2true}
for principally polarised abelian varieties.
The only difference is that for $p=2$ one needs the polarisation to
come from a $K$-rational divisor.
\end{remark}

\begin{remark}
For elliptic curves, the local terms $c_v$ and $|\omega/\neron{v}|_v$
can be obtained from Tate's algorithm, so the conjecture gives an
explicit relation between regulators. Note also that the advantage
of working with regulators up to rational squares is that
one may compute the height matrix on an arbitrary $\Q$-basis
of $E(k)\tensor\Q$.
\end{remark}

\begin{remark}
\label{remCc}
If $E/K$ is semistable, then $C(E/k)$ may be replaced by
just the product of the local Tamagawa numbers $c(E/k)$ in
\ref{C1}-\ref{C2true}.
Indeed, it suffices to show that above a given prime $v$ of $K$, the
contribution from the differential to
$\prod_i C(E/L_i) / \prod_j C(E/L'_j)$ is trivial. But this contribution
is easily seen to be the same for every choice of a local differential
$w_v/K_v$, and it is 1 if $w_v$ is minimal
(as it stays minimal in every extension).
\end{remark}

\begin{remark}\label{finsha}
In \ref{main} and \ref{C2true},
the assumption that $\sha[p^\infty]$ is finite for $A$ over all $L_i,L'_j$
follows from its finiteness over their compositum:
if $A/K$ is an
abelian variety and $L/K$ a finite extension with $\sha(A/L)[p^\infty]$
finite, then $\sha(A/K)[p^\infty]$ is also finite. Indeed, the Weil
restriction of scalars $W_{L/K}(A)$
after an isogeny contains $A$ as a direct summand.
%Since $\sha(A/L)[p^\infty]\iso \sha(W_{L/K}(A)/K)[p^\infty]$
%is assumed finite, so is $\sha(A/K)[p^\infty]$.
Since, by assumption, $\sha(A/L)[p^\infty]\iso \sha(W_{L/K}(A)/K)[p^\infty]$
is finite, so is $\sha(A/K)[p^\infty]$.
\end{remark}

\subsection{Regulator quotients and ranks}
\label{ssmoreregs}

\def\lara{\langle,\rangle}

We now explain how to turn the re\-gulator quotients
from the $\square$-Conjecture into parities of Mordell-Weil ranks.
If $A/K$ is an abelian variety and $\Gal(F/K)\iso G$, consider the decomposition
\hbox{$A(F)\tensor_\Z\Q\iso\oplus\rho_k^{\oplus n_k}$} into $\Q$-irreducible
rational $G$-representations. We will show that for given $L_i,L'_j\subset F$
the regulator quotient equals $\prod_k\RC(\rho_k)^{n_k}$
for purely representation-theoretic quantities $\RC(\rho_k)$
(regulator constants) that do not depend on $A$ or the height pairing.

Let $G$ be a finite group, and $\cH$ a set of
representatives of the subgroups of $G$ up to conjugacy.
Its elements are in one-to-one correspondence with transitive $G$-sets
via $H\mapsto G/H$.
We call an element of $\Z\cH$,
$$
  \rel = \sum\nolimits_i H_i - \sum\nolimits_j H'_j \qquad (H_i, H'_j\in\cH)
$$
a {\em relation between permutation representations\/} if
$\oplus_i\C[G/H_i]\iso\oplus_j\C[G/H'_j]$.
If $\Gal(F/K)\iso G$, then in terms of the fixed fields
$L_i=F^{H_i}$ and $L'_j=F^{H'_j}$,
$$
  \bigoplus_i \Ind_{L_i/K} \triv{L_i} \iso
  \bigoplus_j \Ind_{L'_j/K} \triv{L'_j}.
$$

\begin{notation}
Suppose $V$ is a complex representation of $G$, given with a\linebreak
\hbox{$G$-invariant} non-degenerate Hermitian inner product $\lara$ and a
basis $\{e_i\}$.  We write $\det\lara$
or $\det(\lara|V)$ for the determinant of the matrix
$(\langle e_i, e_j\rangle)_{ij}$.
If $V$ is defined over $\Q$, then the class of $\det\lara$ in $\R^*/\Q^{*2}$
does not depend on the choice of a rational basis.
\end{notation}

\begin{definition}
For each $\Q$-irreducible rational representation $\rho$ of $G$
fix a $G$-invariant real-valued symmetric positive definite
inner product $\lara$ on it, and define the {\em regulator constant\/}
$$
  \RC(\rel,\rho)=\frac{\prod_i \det(\frac{1}{|H_i|}\lara|\rho^{H_i})}
                      {\prod_j \det(\frac{1}{|H'_j|}\lara|\rho^{H'_j})}
    \in\Q^*/\Q^{*2}.
$$
It follows from the theorem below that this is independent of
the choice of the inner product. (In particular, $\RC(\rel,\rho)$ is indeed
in $\Q^*/\Q^{*2}$ rather than $\R^*/\Q^{*2}$, as we can choose
$\lara$ to be $\Q$-valued.)
\end{definition}

\begin{theorem}
\label{regconst}
For any $V\iso \bigoplus_k \rho_k^{n_k}$ with $\rho_k$
rational $\Q$-irreducible representations,
$$
\frac{\prod_i \det(\frac{1}{|H_i|}\lara|V^{H_i})}
     {\prod_j \det(\frac{1}{|H'_j|}\lara|V^{H'_j})} =
\prod_k \RC(\rel,\rho_k)^{n_k} \pmod{\Q^{*2}},
$$
for any $G$-invariant real-valued symmetric positive definite
inner product $\lara$ on $V$.
\end{theorem}

\begin{corollary}
\label{correg}
Let $A/K$ be a principally polarised abelian variety,
and let $\rel$ and $F/K, L_i, L'_j$ be as above.
Let $\{\rho_k\}_k$ be the set of $\Q$-irreducible rational representations of $G$,
and let $n_k$ be the multiplicity of $\rho_k$ in $A(F)\tensor_\Z\Q$. Then
$$
  \frac{\prod_i \Reg(A/L_i)}{\prod_j \Reg(A/L'_j)} =
     \prod_k \RC(\rel,\rho_k)^{n_k} \pmod{\Q^{*2}}.
$$
\end{corollary}

In the remainder of \S\ref{ssmoreregs} we prove Theorem \ref{regconst}.
It suffices to show that the left-hand
side is independent of the choice of an inner product.

\begin{lemma}
\label{detsind}
Let $V$ be a (complex) vector space and $\lara_1, \lara_2$
Hermitian inner products. Then
${\det \lara_1}/{\det \lara_2}$
is independent of the choice of a basis of $V$.
\end{lemma}

\begin{proof}
Changing the basis converts the matrix $X$ of an inner product to $M^t X \bar M$,
where $M$ is the matrix of the basis. The assertion follows from
taking the quotient of the determinants.
\end{proof}

\begin{lemma}
\label{dimsaddup}
Let $\rel=\sum_i H_i - \sum_j H'_j$ be a relation between permutation
representations and $\rho$ a complex representation. Then
$$
  \sum\nolimits_i\dim\rho^{H_i}-\sum\nolimits_j\dim\rho^{H'_j}=0.
$$
\end{lemma}

\begin{proof}
Writing $\lara_G$ for the usual inner product on the space of
characters,
$$
  \textstyle
  \sum\dim\rho^{H_i}=\sum \langle\Res_{H_i}\rho,\triv{H_i}\rangle_{H_i}
  =\sum \langle \rho,\Ind^{G}\triv{H_i} \rangle_G
  = \langle \rho,\oplus\Ind^{G}\triv{H_i} \rangle_G.
$$
There is a similar expression for $H'_j$
and the right-hand sides of the two are the same.
\end{proof}

\begin{lemma}
\label{detrelbasis}
Let $\rel=\sum_i H_i - \sum_j H'_j$ be as above, and
$\tau$ a complex irreducible representation with a
Hermitian $G$-invariant inner product.
For each subgroup $H$ fix a basis of $\tau^H$ and let $M_H$ be the
matrix of the inner product on this basis.
Suppose $\rho\iso\tau^n$ with some Hermitian $G$-invariant inner product
$\lara$.
With respect to the bases of $\rho^H$ induced from
those of $\tau^H$ by this isomorphism,
$$
  \frac{\prod\det (\frac{1}{|H_i|}\lara|\rho^{H_i})}
       {\prod\det (\frac{1}{|H'_j|}\lara|\rho^{H'_j})}=
  \left(\frac{\prod\det \frac{1}{|H_i|} M_{H_i}}
             {\prod\det \frac{1}{|H'_j|} M_{H'_j}}\right)^n \>.
$$
In particular the expression is independent of the choice
of $\lara$.
\end{lemma}

\begin{proof}
Since the Hermitian $G$-invariant inner product on $\tau$ is unique
up to a scalar, the matrix of $\lara$ on $\rho^{H}$ with respect to
the induced basis is
$$
  \begin{pmatrix}
     \lambda_{11} M_H & \lambda_{12} M_H & \ldots & \lambda_{1n} M_H \cr
     \lambda_{21} M_H & \lambda_{22} M_H & \ldots & \lambda_{2n} M_H \cr
     \vdots & \vdots & \ddots & \vdots \cr
     \lambda_{n1} M_H & \lambda_{n2} M_H & \ldots & \lambda_{nn} M_H \cr
  \end{pmatrix} \>,
$$
for some $n\times n$ matrix $\Lambda=(\lambda_{xy})$ not depending on $H$.
Hence
$$
  \det (\tfrac{1}{|H|}\lara|\rho^H) =
    (\det\Lambda)^{\dim \tau^H} (\det \tfrac{1}{|H|} M_H)^{n}\>.
$$
The dimensions $\dim \tau^H$ cancel in $\rel$ by Lemma \ref{dimsaddup},
and the result follows.
\end{proof}

\begin{theorem}
Let $\rel=\sum_i H_i - \sum_j H'_j$ be as above, and
$\rho$ a complex representation of $G$.
Suppose $\lara_1,\lara_2$ are two Hermitian $G$-invariant
inner products on $\rho$.
For each subgroup $H$ fix a basis of $\rho^H$.
Then, computing with respect to these bases,
$$
  \frac{\prod_i\det (\tfrac{1}{|H_i|}\lara_1|\rho^{H_i})}
       {\prod_j\det (\tfrac{1}{|H'_j|}\lara_1|\rho^{H'_j})}=
  \frac{\prod_i\det (\tfrac{1}{|H_i|}\lara_2|\rho^{H_i})}
       {\prod_j\det (\tfrac{1}{|H'_j|}\lara_2|\rho^{H'_j})}\>.
$$
\end{theorem}

\begin{proof}
For each subgroup $H$ and each isotypical component $\rho_l\iso\tau_l^{n_l}$
of $\rho$,
choose a basis of $\tau_l^H$ and induce a basis of $\rho_l^H$ as in the
previous lemma, so
$$
  \frac{\prod_i\det (\tfrac{1}{|H_i|}\lara_1|\rho_l^{H_i})}
       {\prod_j\det (\tfrac{1}{|H'_j|}\lara_1|\rho_l^{H'_j})}=
  \frac{\prod_i\det (\tfrac{1}{|H_i|}\lara_2|\rho_l^{H_i})}
       {\prod_j\det (\tfrac{1}{|H'_j|}\lara_2|\rho_l^{H'_j})}\>.
$$
The isotypical components of $\rho$ are pairwise orthogonal, so taking direct
sums gives the same formula with $\rho$ in place of $\rho_l$.
Finally, applying Lemma \ref{detsind} for every $H_i$, $H'_j$ shows that
we could take any basis of $\rho^{H_i}$, $\rho^{H'_j}$ instead of the
constructed one.
\end{proof}

As a consequence we deduce Theorem \ref{regconst}:
if $\rho$ is rational, and we work up to rational squares,
then we do not have to compute $\det (\tfrac{1}{|H|}\lara_1|\rho^{H})$ and
$\det (\tfrac{1}{|H|}\lara_2|\rho^{H})$ in the same basis.

\subsection{Regulator constants: examples}
\label{ssregex}

\begin{example}
\label{exs3}
The group $G=S_3$ has 3 irreducible representations, namely
$\triv{}$ (trivial), $\epsilon$ (sign) and $\Delta$ (2-dimensional), and
$\cH=\{1,C_2,C_3,S_3\}$.
The submodule of $\Z\cH$ of relations is generated by the following
element $\rel$, with regulator constants
$$
\begin{array}{l|lllll}
& \rlap{$\triv{}$} & \rlap{$\epsilon$} & \rlap{$\Delta$}\cr
\hline
\rel=2S_3+1-2C_2-C_3 & 3 & 3 & 3 \cr
\end{array}
$$
Hence, if $\Gal(F/K)\iso S_3$ and $A/K$ is principally polarised, then
$$
  \Reg(A/K)^2 \Reg(A/F) \Reg(A/F^{C_2})^{-2} \Reg(A/F^{C_3})^{-1}
     = 3^{n_{\triv{}}} 3^{n_\epsilon} 3^{n_\Delta} \cdot \square\>,
$$
with $n_\rho$ the multiplicity of $\rho$ in
$A(F)\tensor_\Z\Q$.
So the parity of ${n_{\triv{}}}+{n_\epsilon}+{n_\Delta}$
(equivalently of $\rkalg AK+\rkalg A{F^{C_3}}+\rkalg A{F^{C_2}}$) is ``computable'',
that is it can be determined from the local invariants
using the $\square$-Conjecture: it is given by
$$
  \ord_3\> C(A/K)^2 C(A/F) C(A/F^{C_2})^{-2} C(A/F^{C_3})^{-1} \mod 2.
$$
This generalises Example \ref{ex2}.
\end{example}

\begin{example}
\label{exa5}
Take $G=A_5$. Here the irreducible rational representations are
${\triv{}},\rho_6,\rho_4,\rho_5$ of dimensions $1,6,4$ and 5, respectively,
and the subgroups are $\cH=\{1,C_2,C_3,C_2\!\times\!C_2,C_5,S_3,D_{10},A_4,A_5\}$.
The lattice of relations is generated by 5 elements, and here are the regulator
constants:
$$
\begin{array}{l|ccccccc}
&{\triv{}} & \rho_6 & \rho_4 & \rho_5\cr
\hline
\rel_1 = 1 - 3\,C_2 + 2\,C_2\!\times\!C_2            &2&1&1&2\cr
\rel_2 = C_2\!\times\!C_2 -2\,D_{10} - A_4 + 2\,A_5  &3&1&3&3\cr
\rel_3 = S_3 -D_{10} - A_4 + A_5                     &3&1&3&3\cr
\rel_4 = 1 - 2\,C_2 - C_5 +2\,D_{10}                 &5&5&5&1\cr
\rel_5 = C_3 - C_5 -2\,A_4 + 2\,A_5                  &15&5&15&3\cr
\end{array}
$$
\end{example}
If $E/K$ is an elliptic curve, it follows that the ``computable''
combinations are $\triv{}+\rho_5, \triv{}+\rho_4+\rho_5$ and
$\triv{}+\rho_6+\rho_4$. (For a general principally polarised
abelian variety only the last two are, see Remark \ref{genC2}.)
For instance, from $\triv{}+\rho_5$, the parity of
$\rkalg E{F^{D_{10}}}$ can be determined from the local invariants.

It is interesting to note that $A_5$ is the only group of order $<120$
for which there is a computable combination of representations
($\triv{}+\rho_6+\rho_4$) where the dimensions add up to an odd number.

\begin{example}
\label{exg20}
Let $G = \smallmatrix 1*0* \subset \GL_2(\F_p)$ for some fixed odd prime $p$.
We write $C_p=\smallmatrix1*01$ and $C_{p-1}=\smallmatrix100*$.
The group $G$ has $p-1$ one-dimensional complex representations
whose direct sum is $\Ind^G\triv{C_{p}}$, and one other $(p\!-\!1)$-dimensional
irreducible representation~$\rho$, namely $(\Ind^G\triv{C_{p-1}})\ominus\triv{G}$.
There is a relation
$$
   \rel = 1 - \>(p\!-\!1)\>C_{p-1} - C_p + (p\!-\!1)\>G.
$$
We have $\RC(\rel,\triv{})=p$ and for $\Q$-irreducible rational
$\sigma\subset(\Ind^G\triv{C_{p}})\ominus\triv{G}$,
$$
  \sigma^G=\sigma^{C_{p-1}}=0, \qquad \sigma^1=\sigma^{C_p}=\sigma,
$$
so $\RC(\rel,\sigma)=p^{\dim\sigma}$. It remains to determine $\RC(\rel,\rho)$. We have
$$
  \rho^G=\rho^{C_{p}}=0, \qquad \rho^1=\rho, \qquad \dim\rho^{C_{p-1}}=1.
$$
If $v\in\rho^{C_{p-1}}$ is non-zero, then $v,gv,\ldots,g^{p-2}v$ is a basis
for $\rho$ with $g=\smallmatrix 1101$.
Note that $v+gv+\ldots+g^{p-1}v=0$ since it is in $\rho^{C_p}$.
Take the Hermitian \linebreak $G$-invariant inner
product on $\rho$ with $\langle v,v\rangle=1$, and let us compute
the matrix $X=(\langle g^n v,g^m v\rangle)$. By $G$-invariance we only need
$\langle v,g^m v\rangle$, but
$$
  0 = \langle v,\sum_{m=0}^{p-1} g^m v\rangle =
    \sum_{m=0}^{p-1} \langle v,g^m v\rangle
$$
and the terms in the right-hand side are equal for $m\neq 0$
by $C_{p-1}$-invariance.
So $\langle v,g^m v\rangle=-\frac{1}{p-1}$, and
$$
  X=
  \begin{pmatrix}
     1 & -\frac{1}{p-1} & \ldots & -\frac{1}{p-1} \cr
     -\frac{1}{p-1} & 1 & \ldots & -\frac{1}{p-1} \cr
     \vdots & \vdots & \ddots & \vdots \cr
     -\frac{1}{p-1} &  -\frac{1}{p-1} & \ldots & 1 \cr
  \end{pmatrix} \>.
$$
The determinant of $X$ is $\frac{p^{p-2}}{(p-1)^{p-1}}$ and it follows that
$\RC(\rel,\rho)=p$. (For $p=3$ this recovers Example \ref{exs3}.)
\end{example}

\begin{corollary}
\label{rkmwft}
Suppose $F/K$ is a Galois extension with Galois group
$\smallmatrix 1*0* \subset \GL_2(\F_p)$,
write $M$ for the fixed field of the commutator subgroup~
$\smallmatrix 1*01$ and $L$ for the fixed field of $\smallmatrix 100*$.
For any principally polarised abelian variety $A/K$ with finite $\sha(A/F)[p^\infty]$,
$$
  \rkalg AK+\rkalg AM+\rkalg AL \equiv \ord_p
  \frac{C(A/F)}{C(A/M)}
  \pmod 2\>.
$$
\end{corollary}

\begin{proof}
Decompose
$A(F)\tensor_\Z\Q=\triv{}^{n_1}\oplus\rho^{n_\rho}\oplus\bigoplus_i\sigma_i^{n_{\sigma_i}}$
into rational irreducibles, with $\sigma_i\subset(\Ind^G\triv{C_{p}})\ominus\triv{G}$.
Combining the above example with Corollary \ref{correg} and
Corollary \ref{C2true} (and Remark \ref{genC2}),
we obtain
$$
  n_1+n_\rho+\sum\nolimits_i n_{\sigma_i}\dim{\sigma_i} \equiv
    \ord_p \frac{C(A/F)C(A/K)^{p-1}}{C(A/M)C(A/L)^{p-1}} \pmod 2.
$$
Finally, $\rkalg AK\!=\!n_1$, $\rkalg AM\!=\!n_1\!+\!\sum n_{\sigma_i}\dim{\sigma_i}$
and \hbox{$\rkalg AL\!=\!n_1\!+\!n_\rho$.}
\end{proof}

\section{Tamagawa numbers and root numbers for elliptic curves}
\label{sroot}

\subsection{Review of root numbers}
\label{ssroot}

We now turn to Tamagawa numbers and their relation to
root numbers, in the special case of elliptic curves.
We refer to \cite{RohV, Kob} for the classification of root numbers
of elliptic curves in odd residue characteristic. Incidentally,
while proving Proposition \ref{thmcv} we came upon the following formula
(case (4)) for local root numbers.
It summarises \cite{Kob}  Thm 1.1 (i), (ii) and Remark 1.2 (ii), (iii).

\newpage

\begin{theorem}
\label{thmrootno}
Let $E/K_v$ be an elliptic curve over a local field. Then
\begin{enumerate}
\item $w(E/K_v)=-1$ if $v|\infty$ or $E$ has split multiplicative reduction.
\item $w(E/K_v)=1$ if $E$ has either good or non-split multiplicative reduction.
\item $w(E/K_v)=(\frac {-1}k)$ if $E$ has additive, potentially multiplicative reduction,
and the residue field $k$ of $K_v$ has characteristic $p\ge 3$.
\item $w(E/K_v) = (-1)^{\lfloor \frac{\ord_v(\Delta)|k|}{12}\rfloor}$, if
$E$ has potentially good reduction,
and the residue field $k$ of $K_v$ has characteristic $p\ge 5$.
Here $\Delta$ is the minimal discriminant of $E$,
and $\lfloor x \rfloor$ is the greatest integer $n\le x$.
\end{enumerate}
\end{theorem}

\begin{proof}
(1,2,3)
This follows from the results of \cite{RohV}, \cite{Kob}.

(4)
Since $p\ge 5$, we have $\ord_v(\Delta)\in\{0,2,3,4,6,8,9,10\}$,
and it specifies the Kodaira-N\'eron reduction type of $E$.
Moreover,
the class of $|k|$~\hbox{modulo}~24 determines the quadratic residue symbols
$(\frac{-1}{k})$, $(\frac{2}{k})$ and $(\frac{3}{k})$. Because in our case $w(E/K_v)$
only depends on the reduction type, $\ord_v(\Delta)$ and these
symbols (\cite{Kob}, Thm 1.1, Remark 1.2), this reduces the proof to a
(short) finite computation.
\end{proof}

\begin{remark}\label{remroot}
In cases (3) and (4) we have the following results,
which are elementary to verify:
\begin{itemize}
\item[(a)] The local root number is unchanged in a totally ramified extension
of degree prime to 12, and
\item[(b)] If the residue field
has square order, then $w(E/K_v)=1$.
\end{itemize}
\end{remark}

\subsection{The case of $\smallmatrix1*0*$-extensions}
\label{sstamagawa}
As in Example \ref{exg20} and Corollary \ref{rkmwft},
suppose $F/K$ has Galois group $G=\smallmatrix1*0*\subset\GL_2(\F_p)$
for some odd prime $p$, and
let $M$ and $L$ be the fixed fields of $\smallmatrix1*01$ and
$\smallmatrix100*$, respectively.
Fix an elliptic curve $E/K$ with an invariant differential $\omega$.
For a prime $v$ of $K$, and $k=K,L,M,F$ set
$$
  W_v(k)=\prod_{\nu|v} w(E/k_\nu), \qquad
  C_v(k)=\prod_{\nu|v} c_\nu \Bigl|\frac{\omega}{\neron{\nu}}\Bigr|_\nu\>,
$$
where $\neron{\nu}$ is a N\'eron differential for $E$ at a prime $\nu$ of $k$.
For $v|\infty$ we define $W_v(k)$ by the same formula and set $C_v(k)=1$.

\begin{proposition}
\label{thmcv}
With fields as above,
let $E/K$ be an elliptic curve with a chosen invariant differential $\omega$,
and let $v$ be a place of $K$.
If $v|6$ and ramifies in $L/K$, assume that $E$ is semistable at $v$.
Then
$$
  \ord_p\frac{C_v(F)\,C_v(K)^{p-1}}{C_v(M)\,C_v(L)^{p-1}}
  \equiv0\pmod2\>\>\>
    \Longleftrightarrow
  \>\>
  W_v(K)W_v(M)W_v(L)=1\>.
$$
\end{proposition}

\begin{proof}
Clearly the left-hand side is the same as
$\ord_p(C_v(F)/C_v(M))\!\mod 2$. Now consider the following cases depending
on the behaviour of $v$ in the (degree $p$ Galois) extension $F/M$.
Note that this extension is ramified if and only if $v$ is ramified in $L/K$.

Case 1: primes above $v$ in $M$ split in $F/M$ (this includes all Archimedean places).
Then $C_v(F)=C_v(M)^p$, so $C_v(F)/C_v(M)$ is a square.
Under the action of the decomposition group $D_v$ at $v$, the
$G$-sets $G/\Gal(F/L)$ and $(G/\Gal(F/M))\smallcoprod(G/G)$ are isomorphic.
So the number of primes above $v$ with a given ramification and
inertial degree is the same in $L$ as in $M$ plus in~$K$. It follows that
the local root numbers cancel, $W_v(L)=W_v(K)W_v(M)$.

Case 2: $F/M$ is inert above $v$. Then $v$ must be
totally split in $M/K$, by the structure of $\Gal(F/K)$. As the number of
primes above $v$ in $M$ is even, $C_v(F)$ and $C_v(M)$ are both squares,
and $W_v(M)=1$.
Since in this case $L_v/K_v$ is Galois of odd degree, $W_v(L)=W_v(K)$ by
Kramer--Tunnell \cite{KT}, proof of Prop. 3.4.

Case 3: $F/M$ is ramified above $v$ and $E$ is semistable at $v$.
The contributions from $\omega$ cancel modulo squares,
and $W_v(K)=W_v(L)$.
If $E$ has split multiplicative reduction over a prime of $M$ above $v$,
this prime contributes $p$ to $C_v(F)/C_v(M)$ and $-1$ to the root number.
If the reduction is either good or non-split, it contributes to neither.

Case 4: $F/M$ is ramified above $v\nmid 6p$, and $E$ has additive reduction
at $v$. Since $v\nmid p$, there is no contribution from $\omega$,
and $v$ is unramified in $M/K$ (again, using the structure of $\Gal(F/K)$
and the fact that totally and tamely ramified Galois extensions of local
fields are abelian.) In particular, $M$ has either even number of primes
above $v$ or they have even residue field extension. In each case $W_v(M)=1$
by Remark \ref{remroot}(b).
It remains to compare $W_v(K), W_v(L)$ and the Tamagawa numbers.

Case 4a: $p\ne 3$. All the Tamagawa numbers are prime to $p$.
Also, because $(p,12)=1$ and $L_v/K_v$ is totally ramified, the root
numbers $w(E/K_v)$ and $w(E/L_v)$ are equal by Remark \ref{remroot}(a).

Case 4b: $p=3$ and $E$ has reduction type $II, II^*, I_0^*, I_n^*$
(resp. $III, III^*$) over $K_v$, and the reduction becomes
$I_0^*, I_n^*$ (resp. $III, III^*$) over $L_v$. By inspection,
the root numbers $w(E/K_v)$ and $w(E/L_v)$ are given by the
same residue symbol (this is also clear from \cite{Kob} 1.1--1.2), so they cancel.
Also the Tamagawa numbers are coprime to 3 (\cite{Sil2} IV.9.4).

Case 4c: $p=3$ and $E$ has reduction $IV, IV^*$ over $K_v$.
The reduction becomes good over $L$, so $W_v(L)=1$ and $C_v(F)=1$.
Over $K_v$ the root number is 1 if and only if
$-3\in K_v^{*2}$ (\cite{Kob} Remark 1.2 (iii)),
that is if~$\mu_3\subset K$. This in turn is equivalent to $v$ being split
in $M/K$ ($K_v$ has a cubic ramified Galois extension if and only if
$\mu_3\subset K_v$.) Equivalently, there are two primes above $v$ in $M$
and the contribution from the Tamagawa numbers is a square.
In the other case, $M/K$ is inert and $C_v(M)=3$
(\cite{Sil2} IV.9.4, Steps 5, 8).

Case 5: $F/M$ is ramified above $v|p$, $p>3$ and $E$ has additive reduction
at $v$. By Remark \ref{remroot}, $w(E/K_v)=w(E/L_v)$, so we need the parity
of $\ord_p(C_v(F)/C_v(M))$ and $W_v(M)$.

Fix a place $w$ over $v$ in $M$.
We can replace $\omega$ by the N\'eron differential of $E/M$ at $w$, as this
changes $C_v(F)/C_v(M)$ by a number of the form $\lambda^p/\lambda$ (which
is a square), and the parity of its $p$-adic valuation remains unchanged.

Case 5a: $E/M_w$ has semistable reduction.
Our minimal model at $w$ stays minimal in any extension, so
there is no contribution from $\omega$.
The result follows as in Case 3.

Case 5b: $E/M_w$ has additive reduction.
The reduction stays additive over $F$ and, since $p>3$, all the Tamagawa
numbers are prime to $p$. If $M$ has either even number of primes above $v$
or the residue fields have even degree over $\F_p$,
then $W_v(M)=1$
(by Remark \ref{remroot}) and it also follows that the contributions from $\omega$ are
squares.

Thus we may assume that $M_w/K_v$ has even ramification degree,
in particular $E$ has potentially good reduction at $v$
(for otherwise it would be multiplicative at $w$).
We may also assume that there is an
odd number of primes over $v$ in $M$, and their residue fields are of odd
degree over $\F_p$. By Theorem \ref{thmrootno} and the fact that $p^2\equiv 1\mod 24$,
$$
  w(M_w) = (-1)^{\lfloor \frac{\ord_v(\Delta)|k|}{12}\rfloor}
         = (-1)^{\lfloor \frac{\ord_v(\Delta)p}{12}\rfloor},
$$
and the right-hand side exactly measures the contribution from
$\bigl|\frac{\omega}{\neron{\nu}}\bigr|_\nu$ for a prime $\nu|w$ of $F$.
The result follows by taking the product over $w|v$.
\end{proof}

Now we reap the harvest:

\begin{theorem}
\label{g20parity}
Let $p$ be an odd prime. As above, let $F/K$ have Galois group
$G = \smallmatrix1*0* \subset \GL_2(\F_p)$, and let $M$ and $L$ be the fixed
fields of $\smallmatrix1*01$ and $\smallmatrix100*$, respectively.
Let $E/K$ be an elliptic curve such that
\begin{enumerate}
\item The $p$-primary component $\sha(E/F)[p^\infty]$ is finite.
\item $E$ is semistable at primes $v|6$ that ramify in $L/K$.
\end{enumerate}
Then
$$
    \rkalg EK\!+\!\rkalg EM\!+\!\rkalg EL\text{ is even}\>\>\Longleftrightarrow\>\>
    w(E/K)w(E/M)w(E/L)=1\>.
$$
\end{theorem}

\begin{proof}
As $\sha(E/F)[p^\infty]$ is assumed to be finite,
by Remark \ref{finsha} the same is true over $K$, $M$ and $L$.
We now apply the theory from \S\ref{ssquareness} to the relation
$$
  \Ind_{F/K} \triv{F} \oplus  (\Ind_{K/K} \triv{K})^{\oplus p-1} \iso
  \Ind_{M/K} \triv{M} \oplus  (\Ind_{L/K} \triv{L})^{\oplus p-1}.
$$
By Corollary \ref{rkmwft}, the sum
$\rkalg EK\!+\!\rkalg EM\!+\!\rkalg EL$ is congruent to
$\ord_p C(E/F)/C(E/M)$ modulo 2. By Proposition \ref{thmcv}, it is even
if and only~if
$$
  \prod_{{\text{$v$ place of $K$}}}\!\!\! w(E/K_v)
  \prod_{{\text{$v$ place of $L$}}}\!\!\! w(E/L_v)
  \prod_{{\text{$v$ place of $M$}}}\!\!\! w(E/M_v)
  =1\>.
$$
\end{proof}

\subsection{Application to the Parity Conjecture}
\label{ssparity}

In \cite{TV-P} Theorem 2 we established the following result:

\begin{theorem}
\label{isogroot}
Suppose that $E/K$ is semistable at primes above $p$ and has a $K$-rational
isogeny of degree $p$. If $p=2$, assume furthermore that $E$ is not supersingular
at primes above $2$. Then
$$
  \rksel EKp \text{ even} \>\>\Longleftrightarrow\>\> w(E/K)=1\>.
$$
\end{theorem}

As an application of Theorem \ref{g20parity}
to $F=K(E[2])$
we can prove a form of the parity conjecture (Conjecture \ref{parityconj})
without the isogeny assumption:

\begin{theorem}[=Theorem \ref{ithmparity}]
\label{parity}
Let $E/K$ be an elliptic curve. Suppose $E$ is semistable at primes dividing
2 and 3 and not supersingular at primes dividing 2. If $\sha(E/K(E[2]))$ has
finite 2- and 3-primary parts, then
$$
  \rkalg EK \text{ even}\>\>\Longleftrightarrow\>\> w(E/K)=1\>.
$$
\end{theorem}

\begin{proof}
Write $F=K(E[2])$, and note that $\Gal(F/K)\subset\GL_2(\F_2)\iso S_3$.
By Remark \ref{finsha},
the 2- and 3-primary parts of $\sha(E/k)$ are finite
for $K\!\subset\! k\!\subset\! F$.

If $E$ has a $K$-rational 2-torsion point, the result follows from
Theorem~\ref{isogroot}. If $F/K$ is cubic, then $\rkalg EK$ and
$\rkalg EF$ have the same parity, and also $w(E/K)=w(E/F)$, so the result
again follows.

We are left with the case when
$\Gal(F/K)\iso S_3\iso \smallmatrix 1*0*\subset \GL_2(\F_3)$.
Let $M$ be the quadratic extension of $K$ in $F$
and $L$ one of the cubic ones. By the above argument, we know that
\beq
  \rkalg EM \text{ even} & \Longleftrightarrow & w(E/M)=1,\cr
  \rkalg EL \text{ even} & \Longleftrightarrow & w(E/L)=1.\cr
\eeq
On the other hand, by Theorem \ref{g20parity} with $p=3$,
$$
    \rkalg EK\!+\!\rkalg EM\!+\!\rkalg EL\text{ is even}\>\>\Longleftrightarrow\>\>
    w(E/K)w(E/M)w(E/L)=1\>.
$$
\end{proof}

\begin{remark}
\label{selparity}
Instead of assuming that $\sha$ is finite in the theorem one may give a
statement about Selmer ranks,
by replacing the use of
Theorem \ref{g20parity} by Corollary \ref{g20selmer} in the proof.
When $\Gal(F/K)\iso S_3$, the parity of
$$
  \rksel EK3+(\rksel EM3-\rksel EM2)+(\rksel EL3-\rksel EL2)
$$
is given by the root number $w(E/K)$,
unconditionally on finiteness of~$\sha$.
In all other cases it is the parity of $\rksel EK2$.

We would also like to remark that if Theorem \ref{isogroot} can be extended
to curves with arbitrary reduction at $v|2$,
and Proposition \ref{thmcv} to extensions where additive primes $v|6$ are
allowed to ramify, the parity conjecture for all elliptic curves over
number fields would follow from finiteness of $\sha$.
\end{remark}

\section{Selmer Groups}\label{sselmer}

Hitherto our main tool was Corollary \ref{C2true}, relating regulators
to Tamagawa numbers assuming that $\sha$ is finite. In \S\ref{ssselgen}
we extend this to an unconditional statement about Selmer ranks.
We get our results (Theorem \ref{thmsel} and Corollary \ref{corsel})
by tweaking Tate--Milne's proof of the isogeny invariance of
the BSD quotient (\cite{MilA}, \S1.7).
The quotient of regulators becomes replaced by a
quantity $Q$ measuring the effect of an isogeny on Selmer groups.
In \S\ref{ssselisog} we review how to construct isogenies between
products of Weil restrictions of an abelian variety, and in \S\ref{ssselmain}
address the question of
turning $Q$ into Selmer ranks (somewhat analogous to turning regulators
into Mordell-Weil ranks.)
This can be done in fair generality
(Theorem \ref{selmain}, Corollary \ref{selav}),
and we illustrate it for
$S_n$-extensions (Ex.~\ref{selsym}),
$\smallmatrix1*0*$-extensions (\S\ref{ssselft}),
and dihedral extensions (\S\ref{ssseldih}).
As a final application, in \S\ref{sspparity} we establish
the $p$-parity conjecture for elliptic curves over $\Q$.

\newpage

\subsection{Invariance of the BSD-quotient for Selmer groups}
\label{ssselgen}

\begin{definition}
For an isogeny $\psi:A\to B$ of abelian varieties over $K$, let
\beq
  Q(\psi) = & |\coker(\psi: A(K)/A(K)_{\tors} \to B(K)/B(K)_{\tors})|\>\>\times\cr
            &  \>\>\times\>\>
           |\ker(\psi: \sha(A)_{\div}\to \sha(B)_{\div})|\>,\cr
\eeq
where $\sha_{\div}$ denotes the divisible part of $\sha$.
\end{definition}

\begin{lemma}
\label{lemQ}
$Q(\psi)$ is finite and satisfies the following properties:
\begin{enumerate}
\item $Q(\psi'\psi)=Q(\psi)Q(\psi')$ if $\psi: A\to B$ and $\psi':B \to C$ are isogenies.
\item $Q(\psi\oplus\psi')=Q(\psi)Q(\psi')$ if $\psi: A\to B$ and $\psi':A' \to B'$ are isogenies.
\item $Q(\psi)=p^{\rksel AKp}$ if $\psi: A\to A$ is multiplication by $p$.
\item If $\deg\psi$ is prime to $p$, then so is $Q(\psi)$.
\end{enumerate}
\end{lemma}

\begin{proof} (2), (3) clear.
(1) follows from the fact that
$\psi: \sha(A)_{\div}\!\to\!\sha(B)_{\div}$ is surjective, and
$\psi: A(K)/A(K)_{\tors} \to B(K)/B(K)_{\tors}$ is injective.
Next, consider the conjugate isogeny $\psi^c: B\to A$,
so that $\psi^c\psi$ is the multiplication by $\deg\psi$ map
on $A$. From (1) and (3), $Q(\psi^c)Q(\psi)$ is finite,
so $Q(\psi)$ is finite. (4) also follows.
\end{proof}

\begin{theorem}
\label{thmsel}
Let $X, Y/K$ be abelian varieties given with exterior
forms $\omegaX, \omegaY$. Suppose $\phi: X\to Y$ is an isogeny and
$\phi^t: Y^t\to X^t$ its dual. Writing $\sha_0(X/K)$ for $\sha(X/K)$
modulo its divisible part and
$$
  \Omega_X=
  \prod\limits_{\vabove{v|\infty}{\text{real}}}
  {\int_{X(K_v)} |\omegaX|}\cdot
  \prod\limits_{\vabove{v|\infty}{\text{complex}}}
  {\int_{X(K_v)}\! \omegaX\!\wedge\omegaXbar}
$$
and similarly for $Y$, we have
\begingroup
$$
  \frac{|Y(K)_{\tors}|}{|X(K)_{\tors}|}
  \frac{|Y^t(K)_{\tors}|}{|X^t(K)_{\tors}|}
  \frac{C(X/K)}{C(Y/K)}
  \frac{\Omega_X}{\Omega_Y}
  \prod_{p|{\deg}\phi}\!
    \frac{|\sha_0(X)[p^\infty]|}{|\sha_0(Y)[p^\infty]|}
     \!=\!
  \frac{Q(\phi^t)}{Q(\phi)}.
$$
\endgroup
\end{theorem}

\begin{proof}
Recall that $\omegaX$ and $\omegaY$ enter into the definition of
$C(X/K)$ and $C(Y/K)$. The left-hand side of the asserted
equation is independent of the choices of $\omegaX$ and $\omegaY$
by the product formula, so choose $\omegaX=\phi^*\omegaY$.
Note also that $\phi$ is an isomorphism
between the $p$-primary parts of $\sha_0(X)$ and $\sha_0(Y)$
for $p\nmid\deg\phi$, so the product involving $\sha$ may be taken over
any sufficiently large set of primes. (In fact, it is simply
$|\sha(X)|/|\sha(Y)|$ if both groups are finite.)

Now we follow closely Tate--Milne's proof in \cite{MilA}, \S1.7.
If $f$ is a homomorphism of abelian groups with finite kernel and cokernel,
write
$$
  z(f) = \frac{|\ker f|}{|\coker f|}.
$$
For $k\supset K$ denote by $\phi(k): X(k)\to Y(k)$ the map induced
by $\phi$ on \hbox{$k$-rational points}, and similarly for $\phi^t$.
For a sufficiently large set of places $S$ of $K$ (\cite{MilA} I.(7.3.1)),
$$
  \prod_{v\in S} z(\phi(K_v)) =
    \frac{z(\phi(K))}{z(\phi^t(K))}\frac{|\sha[\phi^t]|}{|\sha[\phi]|}.
$$
Moreover, $z(\phi(K_v))$ is the contribution from $v$ to ${C(Y/K)}/{C(X/K)}$
for finite places, and the quotient of the corresponding integrals for
infinite places with our choice for $\omegaX, \omegaY$.
(Milne also relates $z(\phi(K))/z(\phi^t(K))$
to the torsion and the regulators and, assuming finiteness of the
Tate-Shafarevich groups, $|\sha[\phi^t]|/|\sha[\phi]|$ to
$|\sha(Y)|/|\sha(X)|$. This gives the usual formula for the
isogeny invariance of the BSD-quotient.)

It remains to show that for every prime $p$,
$$
  \ord_p
  \frac{z(\phi^t(K))}{z(\phi(K))}\frac{|\sha[\phi]|}{|\sha[\phi^t]|}
    =
  \ord_p
  \frac{Q(\phi)}{Q(\phi^t)}
  \frac{|Y^t(K)_{\tors}|}{|X^t(K)_{\tors}|}
  \frac{|Y(K)_{\tors}|}{|X(K)_{\tors}|}
  \frac{|\sha_0(X)[p^\infty]|}{|\sha_0(Y)[p^\infty]|}.
$$
Take an integer $N=p^m$ large enough, so that it kills both the
$p$-power torsion in $X(K)$ and $Y(K)$ and the $p$-parts of
$\sha_0(X)$ and $\sha_0(Y)$.
Applying Lemma \ref{lemQ} (1, 3),
$$
  Q(p^m\phi)=p^{m\rksel XKp}Q(\phi), \qquad
  Q(p^m\phi^t)=p^{m\rksel {Y^t}Kp}Q(\phi^t).
$$
Since $X, Y$ and their duals all have the same $p^\infty$-Selmer rank
(they are all isogenous),
it suffices to verify the claim for $\psi=p^m\phi$. But
\beq
  \ord_p|\sha[\psi]| &=& \ord_p|\sha_0(X)|\cdot
           |\ker(\psi|_{\smallsha(X)_{\div}})|\cr
  \ord_p|\ker(\psi(K))| &=& \ord_p|X(K)_{\tors}|\cr
  \ord_p|\coker(\psi(K))| &=& \ord_p|Y(K)_{\tors}|\cdot
     \smash{\hbox{$
       |\coker(\frac{X(K)}{X(K)_{\tors}}\>
       {\buildrel\psi\over\to}\>\frac{Y(K)}{Y(K)_{\tors}})|,
     $}}
     \cr
\eeq
and similarly for $\psi^t$.
Combining these together yields the assertion.
\end{proof}

\begin{remark}
For $\phi: E\to E'$ a cyclic isogeny of degree $p$, this gives
$$
  \frac{C(E/K)}{C(E'/K)}
  \frac{\Omega_E}{\Omega_{E'}}
     \equiv
  \frac{Q(\phi^t)}{Q(\phi)}
     \equiv
  Q(\phi^t)Q(\phi)
     = Q([p]) = p^{\rksel EKp} \pmod{\Q^{*2}},
$$
which is a formula of Cassels (see Birch \cite{Bir} or Fisher \cite{Fis}).
\end{remark}

\begin{corollary}
\label{corsel}
Let $E/K$ be an elliptic curve with a chosen $K$-differential~$\omega$.
Suppose $L_i/K, L'_j/K$ are finite extensions such that
$$
  X=\prod_i W_{L_i/K}(E), \qquad Y=\prod_j W_{L'_j/K}(E)
$$
are isogenous. If $\phi: X\to Y$ is an isogeny and $\phi^t$ its dual, then
$$
  \frac{\prod_i C(E/L_i)}{\prod_j C(E/L'_j)} \equiv
  Q(\phi^t)Q(\phi) \pmod{\Q^{*2}}.
$$
The same is true if $E$ is replaced by a principally polarised abelian
variety over $K$, possibly up to a factor of 2 if the polarisation is not
induced by a $K$-rational divisor.
\end{corollary}

\begin{proof}
Inducing exterior forms on $W_{L_i/K}(E)$ and $W_{L'_j/K}(E)$ by $\omega$,
we have $\Omega_X=\Omega_Y$. Moreover, $X\iso X^t$, $Y\iso Y^t$
and the $p$-primary parts of $\sha/\sha_{\div}$
have square order by Cassels--Tate pairing.
\end{proof}

\subsection{Isogenies between products of Weil restrictions}
\label{ssselisog}

To make Corollary \ref{corsel} explicit,
recall from Milne's \cite{MilO} \S2 how to construct isogenies
$$
  X=\prod\nolimits_i W_{L_i/K}(A) \overarrow{\phi} \prod\nolimits_j W_{L'_j/K}(A) = Y
$$
for a principally polarised abelian variety $A/K$.
For an extension $L/K$ write $G_L=\Gal(\Kbar/L)$.
Suppose $\oplus_i\Ind_{L_i/K}\triv{L_i}\iso\oplus_j\Ind_{L'_j/K}\triv{L'_j}$,
and consider
$$
  M_X = \oplus_i \Z[G_K/G_{L_i}], \qquad M_Y = \oplus_j \Z[G_K/G_{L'_j}].
$$
These are $G_K$-modules, and satisfy $M_X\tensor\Q\iso M_Y\tensor\Q$.
In general, if $M$ is such a module with a given identification
$M\iso \Z^n$ (as an abelian group), the composition
$$
  s:\>\> G_K \lar \Aut_{\Z}(M) = \Aut(\Z^n) = \GL_n(\Z) \lar \Aut(A^n)
$$
is an element of $H^1(G_K,\Aut_{\Kbar}(A^n))$. It corresponds to
a unique form of $A^n$ over $K$, that is an abelian variety over $K$ such that
$A^n$ is isomorphic to it via an isomorphism $\psi$ defined over $\Kbar$.
(The relation between $\psi$ and $s$ is $s(\sigma)=\psi^{-1}\psi^\sigma$.)
Milne denotes this form $A\tensor M$, and with this notation
$X=A\tensor M_X$ and $Y=A\tensor M_Y$.

Next, a principal polarisation $\lambda: A\to A^t$ induces one on $A^n$.
So we can view $(A\tensor M)^t$ as a form of $A^n$, which is
seen to be the same as \hbox{$A\tensor\Hom(M,\Z)$}. If $M$ is a
permutation module, there is a natural isomorphism $M\iso\Hom(M,\Z)$,
and it induces a principal polarisation on $A\tensor M$.

Now suppose $f: M_X \to M_Y$ is an isogeny of $G_K$-modules
(a $G_K$-invariant injection with finite cokernel),
viewed as an $n\times n$-matrix with integer coefficients. Then
$$
  \phi_f:\>\>\> X \overarrow{\psi_X^{-1}} A^n \overarrow{f} A^n \overarrow{\psi_Y} Y
$$
is an isogeny of degree $|\!\det f|^{2\dim A}$
defined over $K$ (\cite{MilO}, Prop. 6a), with
the dual isogeny
$$
  \phi_f^t:\>\>\> X^t \overlarrow{(\psi_X^{-1})^t} (A^n)^t \overlarrow{f^t} (A^n)^t \overlarrow{\psi_Y^t} Y^t.
$$
With respect to the above polarisations, $f^t$ is the transpose of $f$
(see e.g. \cite{DN} \S1.6, esp. Lemma 3).

To summarise:
the natural identifications $M_X\iso\Z^n$ and $M_Y\iso\Z^n$ induce
principal polarisations on $X$ and $Y$;
an isogeny $f: M_X\to M_Y$ induces an isogeny $\phi_f: X\to Y$
of degree $|\!\det f|^{2\dim A}$;
suppressing the principal polarisations,
$(\phi_f)^t\phi_f = \phi_{f^t f}$ where $f^t$ is the transposed matrix and
$(\phi_f)^t$ is the dual isogeny.
Thus, the right-hand side in Corollary \ref{corsel} for $\phi_f$
becomes $Q$ of an explicit endomorphism $\phi_{f^t f}$ of $X$.

\subsection{Determining $Q$}
\label{ssselmain}

Fix a principally polarised abelian variety $A/K$ and a finite
extension $F/K$ with Galois group $G$.

Let $\{\rho_k\}_k$ be the set of $\Q$-irreducible rational representations of $G$.
For $\rho\in \{\rho_k\}$ we will write $\rkselrep{A}{\rho}{p}$
for the $p^{\infty}$-Selmer rank of $A\tensor\Lambda$, where
$\Z^{\dim\rho}\iso\Lambda\subset\rho$ is any $G$-invariant lattice.
Since the Selmer rank is the same for isogenous abelian varieties,
this is independent of the choice of the lattice. Moreover, for $K\subset L\subset F$
with $\Q[G/\Gal(F/L)]\iso\oplus \rho_k^{n_k}$,
$$
  \textstyle
  \rksel ALp = \sum_k n_k\rkselrep{A}{\rho_k}{p}.
$$
We want to express $Q$ in terms of these Selmer ranks.

\begin{lemma}
\label{lemQV}
Let $V$ be a rational representation of $G$, and $f\in\Aut_G(V)$.
For any $G$-invariant lattice $\Lambda$ with $\Lambda\tensor_\Z\Q=V$ there
is an integer $m\ge 1$ such that $mf$ preserves $\Lambda$, and
$$
  Q(f) := Q(\phi_{mf:\Lambda\to\Lambda})/Q(\phi_{m:\Lambda\to\Lambda}) \in\Q^*
$$
is independent of $\Lambda$ and $m$. It satisfies the following properties:
\begin{enumerate}
\item $Q(f'f)=Q(f')Q(f)$ for $f, f'\in \Aut_G(V)$.
\item $Q(f\oplus f')=Q(f)Q(f')$ for $f\in \Aut_G(V), f'\in \Aut_G(V')$.
\item $Q(f)=\prod_k p^{n_k \rkselrep A{\rho_k}p}$
      if $f: V\to V$ is multiplication by $p$, with $V=\oplus_k\rho_k^{\oplus n_k}$
      the decomposition into rational irreducibles.
\item Suppose $f\in\Aut_G(V)$ has an irreducible minimal
polynomial with \hbox{$p$-adically} integral coefficients
and $p \nmid \det f$. Then $\ord_p Q(f)=0$.
\end{enumerate}
\end{lemma}

\begin{proof}
Independence of $m$ follows from the (obvious) special case $m|m'$.
Similarly we can reduce to the case $\Lambda_1\subset \Lambda_2$ with $m_1=m_2=m$.
Let $\iota: \Lambda_1\to\Lambda_2$ be the inclusion map, and $n\ge 1$ an integer
such that $n\Lambda_2\subset\Lambda_1$.
Then
$$
  n\circ (\phi_{mf:\Lambda_2\to\Lambda_2}) =
    \iota\circ (\phi_{mf:\Lambda_1\to\Lambda_1}) \circ (n\iota^{-1})\>.
$$
The independence of $\Lambda$
now follows from Lemma \ref{lemQ} (1).

(1-3) are immediate from Lemma \ref{lemQ}.

(4) Let $m(x)$ be the minimal polynomial of $f$.
After scaling $f$ by an integer coprime to $p$ if necessary,
we may assume that $m(x)$ has integer coefficients. Let $\alpha\in\bar{\Q}$
be a root of $m$, and let $\K=\Q(\alpha)$. Note that $\alpha$ is an algebraic
integer, $\alpha\in\O_\K$.

Via the action of $f$ and of $G$, the representation $\rho$
is naturally a $\K[G]$-module. Take a $G$-invariant full $\O_\K$-lattice $\Lambda$.
(It exists, since one may take any full $\O_\K$-lattice,
and generate a lattice by its $G$-conjugates.)
In particular, it is a full $G$-invariant $\Z$-lattice preserved by $f$,
so $Q(f)=Q(\phi_f)/Q(\phi_{\id})$ is an integer.
By Lemma \ref{lemQ}(4) it is prime to $p$, as $p$ does not divide
$\deg\phi_f=|\!\det f|^{2\dim A}$.
\end{proof}

\begin{theorem}
\label{selmain}
Let $V\iso\rho^{\oplus n}$, with $\rho$ a $\Q$-irreducible rational representation
of $G$, and let $f\in\Aut_G(V)$. Suppose $p$ is a prime such that either
\begin{enumerate}
\item\label{condqp} $\rho$ is irreducible as a $\Q_p[G]$-representation, or
\item\label{condmx} for every irreducible factor $m(x)|\det(f-x I)\in \Q[x]$,
all of the roots of $m(x)$ in $\Qpb$ have the same valuation.
\end{enumerate}
Then
$$
  \ord_p Q(f) \equiv
     \frac{\ord_p \det f}{\dim\rho}
     \>
     \rkselrep{A}{\rho}{p}
     \mod 2.
$$
\end{theorem}

\begin{proof}
\eqref{condmx}
First, we can break up $V$ as follows. Let
$$
  \det(f-x I) = \prod_k m_{k}(x)^{n_{k}}
$$
be the factorisation into irreducibles. Then $V_{k}=\ker m_{k}(f)^{n_{k}}$
is $G$-invariant because
$f$ commutes with the action of $G$ (so $G$ preserves its generalised
eigen\-spaces). Since $V\!=\!\oplus V_k$, by Lemma \ref{lemQV} (2)
it suffices to prove
the state\-ment with $V\!=\!V_{k}$ and $\det(f\!-\!x I)$ a power of an irreducible
polynomial~$m(x)$.

Suppose all the roots of $m(x)$ in $\Qpb$ have the same valuation.
Write
$$
  f^{\dim V}=p^{\ord_p\det f}\cdot f',
$$
so that the roots of the characteristic polynomial of $f'$ in $\Qpb$
are all $p$-adic units. Then $\ord_p Q(f')=0$ by Lemma \ref{lemQV} (4), and
the claim follows by Lemma \ref{lemQV} (1, 3).

\eqref{condqp}
Fix an identification $V=\rho^{\oplus n}$.
As $D=\End_G(\rho)$ is a skew field, we can put $f$ into a block-diagonal form
by multiplying it on the left and on the right by $n\times n$
matrices with values in $D$ that are (a) permutation matrices and
(b) identity plus some element of $D$ in $(i,j)$-th place ($i\ne j$).
(This is just the usual Gaussian elimination over a skew field.)
All of these elementary matrices have $Q=1$
(they are either of finite order or commutators)
% for commutators: use \smallmatrix 1101 and \smallmatrix 100e
and $\det=\pm 1$,
so we are reduced to the case $V=\rho$.

We claim that for $V$ irreducible over $\Q_p$,
the eigenvalues of $f$ in $\Qpb$ have the
same valuation (so \eqref{condmx} applies). But otherwise the minimal
polynomial of $f$ is reducible over $\Q_p$,
and we can decompose $V$ over $\Q_p$ as above,
contradicting the irreducibility.
\end{proof}

\begin{corollary}
\label{selav}
Let $K\subset L_i, L'_j\subset F$ be finite extensions with $F/K$ Galois
with Galois group $G$.
Suppose there is an isogeny of $\Z[G]$-modules
$$
  f:\>\> \prod_i \Z[G/H_i] \lar \prod_j \Z[G/H'_j],
$$
where $H_i=\Gal(F/L_i)$ and $H'_j=\Gal(F/L'_j)$. Assume furthermore that
on every isotypical component $\rho^{n_\rho}$ of $\prod_i \Q[G/H_i]$ the automorphism
$f^tf$ satisfies either (1) or (2) of Theorem \ref{selmain}.
Then for every elliptic curve $E/K$ with a chosen $K$-differential~$\omega$,
$$
  \ord_p
  \frac{\prod_i C(E/L_i)}{\prod_j C(E/L'_j)} \equiv
  \sum_\rho \>\>
     \frac{\ord_p \deton{f^tf}{\rho^{n_\rho}}}{\dim\rho}
  \>
     \rkselrep{E}{\rho}{p}\>
  \mod 2,
$$
the sum taken over the distinct $\Q$-irreducible rational representations of $G$.
\end{corollary}

\begin{remark*}
This also holds for principally polarised abelian
varieties for odd~$p$, and for $p=2$ provided the polarisation is
induced by a $K$-rational divisor.
\end{remark*}

Here are some special cases when the theorem applies:

\begin{example}
\label{selsym}
If $G=S_n$, then every $\Q$-irreducible rational representation is absolutely
irreducible, so the condition \eqref{condqp} of Theorem \ref{selmain} always
holds. Thus the corollary applies for every isogeny and all $p$.
\end{example}

\begin{example}
\label{sel60}
For all groups with $|G|\le 55$, every relation of permutation
representations and every prime $p$,
we have checked that there is always an isogeny satisfying
one of the conditions of Theorem \ref{selmain} on every isotypical
component. In each case, the coefficient of $\rkselrep A{\rho}p$ agrees
with the regulator constant of \S\ref{ssmoreregs}. Is this true in general?%
$\,$\footnote{Yes, see \cite{Selfduality} which provides an analogue of
              regulator constants for Selmer groups.}
\end{example}

\subsection{Example: Selmer ranks for $\smallmatrix1*0*$-extensions}
\label{ssselft}

As an illustration, we extend Corollary \ref{rkmwft} to Selmer ranks:

\begin{theorem}
\label{thmselft}
Let $p$ be an odd prime. Suppose $F/K$ has Galois group
$G = \smallmatrix1*0* \subset \GL_2(\F_p)$, and let $M$ and $L$ be the fixed
fields of $\smallmatrix1*01$ and $\smallmatrix100*$, respectively.
For every principally polarised abelian variety $A/K$,
$$
  \rksel AKp+\rksel ALp+\rksel AMp
  \equiv \ord_p\frac{C(A/F)}{C(A/M)} \pmod2.
$$
\end{theorem}

\begin{proof}
Consider the abelian varieties
$$
  X=W_{L/K}(A)^{p-1}\times W_{M/K}(A), \qquad Y=A^{p-1}\times W_{F/K}(A).
$$
By Corollary \ref{corsel}, it suffices to show that
$$
  \rksel AKp+\rksel ALp+\rksel AMp   \equiv  \ord_p Q(\phi^t\phi) \pmod2
$$
for some isogeny $\phi: X\to Y$.
Write $\Gal(F/M)=\langle g\rangle, \Gal(F/L)=\langle h\rangle$ with
$g^p=1=h^{p-1}$, and introduce permutation modules
$$
  \begin{array}{llllll}
  \ZG{K}&=&\Z[G/G]&=&\Z\cr
  \ZG{L}&=&\Z[G/\langle h\rangle]&=&\oplus\Z g^i & \qquad 0\le i\le p-1\cr
  \ZG{M}&=&\Z[G/\langle g\rangle]&=&\oplus\Z h^j & \qquad 0\le j\le p-2\cr
  \ZG{F}&=&\Z[G]&=&\oplus\Z g^ih^j. \cr
  \end{array}
$$
Consider
\beq
  V_1&=&\ZG{L}x_1\oplus\ldots\oplus \ZG{L}x_{p-1}
    \>\>\oplus\>\> \ZG{M}x_p,\cr
  V_2&=&\ZG{K}y_1\oplus\ldots\oplus \ZG{K}y_{p-1}
    \>\>\oplus\>\> \ZG{F}y_p,
\eeq
and take the $G$-invariant map $f: V_1\to V_2$ determined by
\beq
  x_1     &\mapsto& y_1 &+& {\textstyle\sum_j}\, h^j\,y_p \cr
  x_k     &\mapsto& y_1-y_k &+& {\textstyle\sum_j}\, h^j\,(1\!-\!g^{1-k})\> y_p  \quad (k=2,\ldots,p-1) \cr
  x_p     &\mapsto& y_1+\ldots+y_{p-1} &-& {\textstyle\sum_i}\> h^{-1}g^i\> y_p. \cr
\eeq
It is easy to check that it is well-defined, and moreover,
written as a matrix on the chosen $\Z$-basis of $V_1$ and $V_2$,
$$
  |\det f| = (p^2-p+1) p^{\frac{p(p-1)}2-1},
$$
in particular non-zero. So $f$ induces an isogeny $\phi_f: X\to Y$
(\S\ref{ssselisog}). Next, $\phi_f^t$ is given by the transposed matrix
(\S\ref{ssselisog} again),
%\beq
%  y_1     &\mapsto& {\textstyle\sum_i \sum_j}\, g^i\, x_{j+1}
%          &+& {\textstyle\sum_j} h^j x_p\cr
%  y_k     &\mapsto& -{\textstyle\sum_i}\, g^i\,x_k &+& {\textstyle\sum_j} h^j x_p  \quad (k=2,\ldots,p-1)\cr
%  y_p     &\mapsto&  {\textstyle\sum_{i\ne 0}}\,(1\!-\!g^{i-1})\, x_i &-& h^{-1} x_p  \cr
%\eeq
and the composition $f^tf$~by
\beq
  x_1     &\mapsto& {\textstyle\sum_{i\ne 0}}\, (g^i x_1  + p x_i) \cr
  x_k     &\mapsto& px_k + \sum_{i\ne 0}\, px_i&(k=2,\ldots,p-1) \cr
  x_p     &\mapsto& (p+{\textstyle\sum_j} (p-1)h^j)\,x_p\>.
\eeq
(As an example, for $p=3$ the maps $f$ and $f^tf$ are
\begingroup
$$
\font\txtfnt=cmr8\textfont0=\txtfnt
\baselineskip 8pt
\def\arraystretch{0.8}
\def\arraycolsep{5pt}
f:
\left(\!
\begin{array}{ccc|ccc|cc}
1&1&1&1&1&1&1&1\cr
\hline
0&0&0&\llap{-}1&\llap{-}1&\llap{-}1&1&1\cr
\hline
1&0&0&1&\llap{-}1&0&0&\llap{-}1\cr
0&1&0&0&1&\llap{-}1&0&\llap{-}1\cr
0&0&1&\llap{-}1&0&1&0&\llap{-}1\cr
1&0&0&1&0&\llap{-}1&\llap{-}1&0\cr
0&0&1&0&\llap{-}1&1&\llap{-}1&0\cr
0&1&0&\llap{-}1&1&0&\llap{-}1&0\cr
\end{array}
\!\right)
\quad
f^tf:
\left(\!
\begin{array}{ccc|ccc|cc}
3&1&1&3&0&0&0&0\cr
1&3&1&0&3&0&0&0\cr
1&1&3&0&0&3&0&0\cr
\hline
3&0&0&6&0&0&0&0\cr
0&3&0&0&6&0&0&0\cr
0&0&3&0&0&6&0&0\cr
\hline
0&0&0&0&0&0&5&2\cr
0&0&0&0&0&0&2&5\cr
\end{array}
\!\right)
$$
\endgroup
with respect to the bases $\{x_1,gx_1,g^2x_1,x_2,gx_2,g^2x_2,x_3,hx_3\}$
of $V_1$ and $\{y_1,y_2,y_3,gy_3,g^2y_3,hy_3,hgy_3,hg^2y_3\}$
of $V_2$.)

Clearly $f^tf=\alpha_1\oplus\alpha_2$, with
$\alpha_1$ an endomorphism of $\Z_L^{p-1}$ and
$\alpha_2$ of $\Z_M$. To prove the theorem it suffices to show that
\beq
  \ord_p Q(\alpha_1) = {(p-2)\rksel{W_{L/K}(A)}Kp}, \cr
  \ord_p Q(\alpha_2) = {\rksel{W_{M/K}(A)}Kp - \rksel AKp}.
\eeq
The map $\alpha_1$ is
the composition of \hbox{$\id\oplus p\oplus\cdots\oplus p$} with an endomorphism
of determinant $p^2-p+1$. Each multiplication by $p$ on a copy of
$W_{L/K}(A)$ contributes $p^{\rksel{W_{L/K}(A)}Kp}$ to $Q(\alpha_1)$,
and the remaining endomorphism contributes nothing (Lemma \ref{lemQ}).

As for $\alpha_2$, consider
$$
  \alpha_2\oplus [p]:\>\, \ZG{M} z_1  \oplus \ZG{K} z_2 \lar \ZG{M} z_1  \oplus \ZG{K} z_2.
$$
It is easy to check that
$$
  \alpha_3 \>\circ\> (\alpha_2\oplus [p])
    = ([p]\oplus\id) \>\circ\> \alpha_4,
$$
with
$$
  \alpha_3: \left\{ \begin{array}{lll}
    z_1\to z_1+\sum_j h^jz_1 \cr
    z_2\to (p-1)\sum_j h^jz_2\cr
  \end{array}\right., \quad
  \alpha_4: \left\{ \begin{array}{lll}
    z_1\to z_1+p\sum_j h^jz_1+z_2\cr
    z_2\to (p-1)(p^2-p+1)\sum_j h^jz_1\cr
  \end{array}\right..
$$
Furthermore, $\det\alpha_3$ and $\det\alpha_4$ are prime to $p$, and it
follows that $\ord_p Q(\alpha_2)$ is ${\rksel{W_{M/K}(A)}Kp - \rksel AKp}$,
as asserted.
\end{proof}

Using the result on the local Tamagawa numbers
in this extension (Proposition \ref{thmcv}) we
obtain the following strengthening of Theorem \ref{g20parity}.

\begin{corollary}\label{g20selmer}
Let $p$ be an odd prime. As above, let $F/K$ have Galois group
$G = \smallmatrix1*0* \subset \GL_2(\F_p)$, and let $M$ and $L$ be the fixed
fields of $\smallmatrix1*01$ and $\smallmatrix100*$, respectively.
For every elliptic curve $E/K$ with semistable reduction at the primes $v|6$
that ramify in $L/K$,
$$
    \rksel EKp\!+\!\rksel EMp\!+\!\rksel ELp\text{ is even}\>\>\Leftrightarrow\>\>
    w(E/K)w(E/M)w(E/L)=1\>.
$$
\end{corollary}

This
can be used to study the ranks of elliptic curves in an infinite
``false Tate curve extension'' with Galois group
$\smallmatrix1*0* \subset \GL_2(\Zp)$. Arithmetic of elliptic curves
(ordinary at $p$) in such extensions has been studied in the context of
non-commutative Iwasawa theory, see e.g. \cite{HV,CFKS,CS}.

Thus, fix a number field $K$, an odd prime $p$,
and $\alpha\in K^*$. We are interested in the extensions
$K(\sqrt[p^n]{\alpha})$ and $K(\mu_{p^n},\sqrt[p^n]{\alpha})$ of $K$.
We will assume that their degree is maximal possible,
i.e. $p^n$ and $(p-1)p^{2n-1}$, respectively.

\begin{proposition}
\label{ftate}
Let $E/K$ be an elliptic curve for which
the parity of the $p^\infty$-Selmer rank agrees with the root number over
$K$ and over $K(\mu_p)$, and semistable at
all primes $v|6$ that ramify in $K(\sqrt[p^n]{\alpha})/K$.
Then
$$
\rksel E{{K(\sqrt[p^n]{\alpha})}}p
\text{ is even}\>\>\Longleftrightarrow\>\>
w(E/K(\sqrt[p^n]{\alpha}))=1\>.
$$
\end{proposition}

\begin{proof}
For brevity, let us write $L_i=K(\sqrt[p^i]{\alpha})$ and
$F_i=K(\mu_p,\sqrt[p^i]{\alpha})$ for $i\ge 0$.
We prove the result by induction on $i$, by showing that if the parity of
the $p^\infty$-Selmer rank agrees with the root number over $L_{i-1}$ and
$F_{i-1}$ then it does so over $L_i$ and $F_i$ (for $1\le i\le n$).

The extension $F_i/F_{i-1}$ is Galois of odd degree, so
$w(E/F_i)=w(E/F_{i-1})$.
By Corollary \ref{corcycext} below, the parity of the $p^\infty$-Selmer rank
is also unchanged in this extension.
The fact that the parity of the $p^\infty$-Selmer rank agrees with the
root number over $L_i$ follows from Corollary~\ref{g20selmer}
applied to the extension~$F_i/L_{i-1}$.
\end{proof}

The following is a standard result on the behaviour of Selmer groups
in Galois extensions. We give a brief proof for lack of a reference.
Write $\Sel_{p^n}$ and $\Sel_{p^\infty}$ for the $p^n$- and
$p^{\infty}$-Selmer groups, and set
$$
  \X_p(E/K)=\Hom(\Sel_{p^\infty}(E/K),\Qp/\Zp)\tensor_{\Zp}\Qp\>.
$$
The $p^\infty$-Selmer rank of $E/K$ is the same as the dimension of $\X_p(E/K)$
as a $\Qp$-vector space.

\begin{lemma}\label{seletale}
Let $E/K$ be an elliptic curve, and let $F/K$ be a finite Galois extension
with Galois group $G$. Then
$$
  \rksel EKp = \dim_{\Q_p} \X_p(E/F)^G\>.
$$
\end{lemma}

\begin{proof}
The restriction map from $H^1(K,E[p^n])$ to $H^1(F,E[p^n])^G$ induces
a map $\Sel_{p^n}(E/K)\to\Sel_{p^n}(E/F)^G$ whose kernel and cokernel
are killed by $|G|^2$. Taking direct limits gives a map from
$\Sel_{p^\infty}(E/K)$ to $\Sel_{p^\infty}(E/F)^G$,
whose kernel and cokernel are killed by $|G|^2$.
The result follows by taking duals and tensoring with $\Qp$.
\end{proof}

A cyclic group of order $p$ has only two $\Q_p$-irreducible $p$-adic
representations, the trivial one and one of dimension $p-1$. Thus,

\begin{corollary}
\label{corcycext}
The parity of the $p^\infty$-Selmer rank in unchanged in cyclic
$p$-extensions.
\end{corollary}

\begin{example}
\label{ex41a1}
Let $p=3$ and consider the elliptic curve
$$
  E: y^2+xy=x^3-x^2-2x-1\qquad(49\text{A}1).
$$
It has additive reduction of Kodaira type III at 7 and is supersingular at 3.

For a false Tate curve extension, we take
$\Q(\mu_{3^n},\sqrt[3^n]{m})$ for some cube free $m>1$. Using
3-descent for $E$ and its quadratic twist by $-3$ over $\Q$, it is easy
to see that $\rkalg E\Q\!=\!\rksel E\Q3\!=\!0$ and
$\rkalg E{\Q(\mu_3)}\!=\!\rksel E{\Q(\mu_3)}3\!=\!1$, both in agreement
with the root numbers.

By Proposition \ref{ftate}, the $3^\infty$-Selmer rank of $E$ over $L_n=K(\sqrt[p^n]{m})$
agrees with the root number, which equals
$(-1)^n$ for every $m$ ($-1$ from $v|7$ and $(-1)^{n-1}$ from $v|\infty$).
Because the Selmer rank is non-decreasing in extensions (e.g. by Lemma
\ref{seletale}), the $3^\infty$-Selmer rank must be at least
$n$ over~$L_n$.
In fact, using Lemma \ref{seletale}
and that $\Ind_{L_n/\Q}\triv{L_n}\ominus\Ind_{L_{n-1}/\Q}\triv{L_{n-1}}$
is irreducible, it is easy to see that the $3^\infty$-Selmer rank over
$\Q(\mu_{3^n},\sqrt[3^n]{m})$
is at least $3^n$.

%(By Theorem \ref{isogroot} below, for this curve
%the parity of the $2$-Selmer rank also
%agrees with the root number over any field, since $E$ has a 2-torsion
%point $(2,-1)$ over $\Q$ and has good ordinary reduction at 2.)

\end{example}

\subsection{Example: Dihedral groups}
\label{ssseldih}
As another illustration, we consider dihedral groups to obtain similar results
to \cite{MR}, e.g. Theorem 8.5.
For simplicity, we will only look at $D_{2p}$ with $p$ an odd prime.

\begin{proposition}
\label{propdih}
Suppose $\Gal(F/K)=D_{2p}$ with $p$ an odd prime, and pick
extensions $M/K$ and $L/K$ in $F$ of degree 2 and $p$, respectively.
For every principally polarised abelian variety $A/K$,
$$
  \rksel A{M}p+\tfrac{2}{p-1}(\rksel A{L}p-\rksel A{K}p)
  \equiv
  \ord_p\frac{C(A/F)}{C(A/M)}
  \pmod2.
$$
\end{proposition}

\begin{proof}
First let $G=D_{2n}=\langle g,h\>|\>g^n=h^2=hghg=1\rangle$
for a general $n$, and write
$n=2m+\delta$ with $\delta\in\{0,1\}$. Take the permutation modules
\beq
  V_1=v_1\Z[G/\langle g^{-1}h\rangle]\oplus v_2\Z[G/\langle g^{-2}h\rangle] \oplus
        v_3\Z[G/\langle g\rangle],\cr
  V_2=w_1\Z[G/G]\oplus w_2\Z[G/G]\oplus w_3\Z[G].
\eeq
Consider
$f: V_1\to V_2$ and $f^tf: V_1\to V_1$ given respectively by
\beq
  v_1 \mapsto (1+g^{-1}h)w_3 &&& v_1 \mapsto 2v_1\cr
  v_2 \mapsto w_2\!+\!g^{m-2}(1\!-\!g^{1+\delta})(g\!-\!h)w_3
    &&& v_2 \mapsto \bigl(4\!-\!2g^{1+\delta}\!-\!2g^{-1-\delta}
        \!+\!\sum\limits_{i=0}^{\smash{n-1}} g^i\bigr)v_2 \cr
  v_3 \mapsto w_1 + \sum\limits_{i=0}^{\smash{n-1}} g^i (1-h) w_3
    &&& v_3 \mapsto \bigl((2n+1)-(2n-1)h\bigr)v_3 \cr
\eeq
Note that $f^tf$ decomposes naturally as $\alpha_1\oplus\alpha_2\oplus\alpha_3$
on $V_1$, and that $\alpha_2$ is given on the basis
$\{v_2,gv_2,\ldots g^{n-1}v_2\}$ by the matrix
$$
\begingroup
\font\txtfnt=cmr8\textfont0=\txtfnt
\baselineskip 8pt
\def\arraystretch{0.8}
\def\arraycolsep{5pt}
\let\v\vdots
\let\h\hdots
\let\d\ddots
\def\m{\llap{-}1}
\begin{array}{cc}
\left(\!
\begin{array}{ccccccccc}
  5 &  \m &   1 &   1 &   1 &  \h &   1 &   1 &  \m   \cr
 \m &   5 &  \m &   1 &   1 &  \h &   1 &   1 &   1   \cr
  1 &  \m &   5 &  \m &   1 &  \h &   1 &   1 &   1   \cr
  1 &   1 &  \m &   5 &  \m &  \h &   1 &   1 &   1   \cr
  1 &   1 &   1 &  \m &   5 &  \h &   1 &   1 &   1   \cr
 \v &  \v &  \v &  \v &  \v &  \d &  \v &  \v &  \v   \cr
  1 &   1 &   1 &   1 &   1 &  \h &   5 &  \m &   1   \cr
  1 &   1 &   1 &   1 &   1 &  \h &  \m &   5 &  \m   \cr
 \m &   1 &   1 &   1 &   1 &  \h &   1 &  \m &   5   \cr
\end{array}
\!\right)
  &
\left(\!
\begin{array}{ccccccccc}
  5 &   1 &  \m &   1 &   1 &  \h &   1 &  \m &   1   \cr
  1 &   5 &   1 &  \m &   1 &  \h &   1 &   1 &  \m   \cr
 \m &   1 &   5 &   1 &  \m &  \h &   1 &   1 &   1   \cr
  1 &  \m &   1 &   5 &   1 &  \h &   1 &   1 &   1   \cr
  1 &   1 &  \m &   1 &   5 &  \h &   1 &   1 &   1   \cr
 \v &  \v &  \v & \v  &  \v &  \d &  \v &  \v &  \v   \cr
  1 &   1 &   1 &   1 &   1 &  \h &   5 &   1 &  \m   \cr
 \m &   1 &   1 &   1 &   1 &  \h &   1 &   5 &   1   \cr
  1 &  \m &   1 &   1 &   1 &  \h &  \m &   1 &   5   \cr
\end{array}
\!\right)\cr
\vphantom{X^{{\int}^X}}
\text{$n\ge\> $}4\>\>\text{even} &
\text{$n\ge\> $}3\>\>\text{odd}
\end{array}
\endgroup
$$
with $\det\alpha_2=2^{n-1}n^3$ for any $n\ge 2$.

Now suppose that $n=p$ is an odd prime, so
$$
  A\tensor V_1=W_{L/K}(A)^2 \times W_{M/K}(A), \qquad
  A\tensor V_2=A^2 \times W_{F/K}(A).
$$
By Corollary \ref{corsel}, it suffices to show that $\ord_p Q(\phi_{f^tf})$
has the same parity as the left-hand side of the formula in the
proposition. Since $f^tf=\alpha_1\oplus\alpha_2\oplus\alpha_3$,
it remains to determine $\ord_p Q(\alpha_i)$.
Clearly $Q(\alpha_1)$ is prime to $p$. Next,
$\alpha_3$ acts as multiplication by $2$ (resp. $4p$) on the trivial
(resp. ``sign'') component of $\Q[G/\langle g\rangle]$, so
$\ord_p Q(\alpha_3)=\rksel A{M}p-\rksel A{K}p$.

Finally, $\alpha_2$ on $\Q[G/\langle g^{-2}h\rangle]\iso\triv{}\oplus\rho$
has determinant $p$ on $\triv{}$ and therefore
determinant $2^{p-1}p^2$ on $\rho$.
As $\triv{},\rho$ are $\Q_p$-irreducible, Theorem \ref{selmain}
applies:
$$
  \ord_p Q(\alpha_2) = \rkselrep A{\triv{}}p + \tfrac{2}{p-1}
    \rkselrep A{\rho}p\>.
$$
Since $\rkselrep A{\rho}p = \rksel ALp - \rksel AKp$, this
completes the proof.
\end{proof}

\begin{remark}
\label{dihsimple}
Let $E/K$ be an elliptic curve, and for simplicity let $p>3$. Then
$$
  \ord_p \frac{C(E/F)}{C(E/M)} \equiv
    |S_1|+|S_2| \mod 2,
$$
where $S_1$ (resp. $S_2$) is the set of primes $v$ of $M$ that ramify in $F/M$
where $E$ has split multiplicative reduction (resp. additive reduction,
$v|p$, $M_v/\Q_p$ has odd residue degree,
and $\lfloor p\ord_v(\Delta_v)/12\rfloor$ is odd; $\Delta_v$ is the
minimal discriminant of $E$ at $v$).
So if $\rksel EMp+|S_1|+|S_2|$ is odd, then
$$
  \rksel ELp \ge \rksel EKp + \tfrac{p-1}2.
$$
\end{remark}

\subsection{Application to the $p$-Parity Conjecture over $\Q$}
\label{sspparity}

\begin{theorem}[=Theorem \ref{ithmpparity}]
\label{pparity}
For every elliptic curve $E/\Q$ and every prime $p$,
$$
  \rksel E\Q p \equiv \ord_{s=1}L(E,s)\mod 2\>.
$$
\end{theorem}

\begin{proof}
For $p=2$ this is due to Monsky \cite{Mon}, so suppose $p$ is odd.
(Presumably the proof below would work for
modular abelian varieties over totally real fields.)

By the results of Bump--Friedberg--Hoffstein--Murty--Murty--Waldspurger\\
\cite{BFH,MM,Wal},
there is an imaginary quadratic field $M_0$ where all bad primes of $E$
split, and such that the quadratic twist of $E$ by $M_0$
has analytic rank at most 1. By Kolyvagin's theorem \cite{Kol},
the parity conjecture holds for the twist, so it suffices to prove
it for $E/M_0$.

Let $M_n$ denote the $n$-th layer in the anticyclotomic $\Z_p$-extension
of $M_0$. The parity of the analytic rank is the same over $M_n$ as over
$M_0$ since the root number is unchanged in cyclic odd-degree extensions.
The same holds for the $p^{\infty}$-Selmer rank (Corollary \ref{corcycext}),
so it suffices to prove the parity statement for $E/M_n$ for some $n$.

Embedding $M_{n+1}$ in $\C$, complex conjugation acts
on the cyclic group $\Gal(M_{n+1}/M_0)$ as $-1$, so $\Gal(M_{n+1}/\Q)$
is dihedral. Write
$$
  F=M_{n+1}, \quad
  M=M_n, \quad
  L=M_{n+1}\cap\R, \quad
  K=M_{n}\cap\R.
$$
Then $H=\Gal(F/K)\iso D_{2p}$. It has three $\Q_p$-irreducible $p$-adic
representations: trivial $\triv{}$, sign $\epsilon$ and
$(p\!-\!1)$-dimensional $\rho$.
As before, write $\X_p(E/k)$ for
the dual Selmer group
$\Hom(\Sel_{p^\infty}(E/k),\Qp/\Zp)\tensor_{\Zp}\Qp$, and
decompose
$$
  \X=\X_p(E/F)\>\iso\>
    \triv{}^{\oplus m_{\triv{}}}\oplus
    \epsilon^{\oplus m_\epsilon}\oplus
    \rho^{\oplus m_\rho}.
$$
As $\X_p(E/K)=\X^H$ etc. (Lemma \ref{seletale}),
$$
  \rksel EKp=m_{\triv{}}, \quad
  \rksel EMp=m_{\triv{}}\!+\!m_\epsilon, \quad
  \rksel ELp=m_{\triv{}}\!+\!\tfrac{p-1}2 m_\rho.
$$
Now we invoke Proposition \ref{propdih}:
$$
  \rksel EMp + {m_\rho} \equiv \ord_p\tfrac{C(E/F)}{C(E/M)} \mod 2.
$$
Since all bad primes of $E$ split in $M/K$,
the root number $w(E/M)=-1$ and both $C(E/F)$ and $C(E/M)$ are squares.
So the right-hand side in the above formula is zero, and
it suffices to show that $m_\rho$ is odd.

Now take $n$ large enough. Then Cornut--Vatsal's \cite{CV} Thm. 1.5
provides a primitive character~$\chi$
of $\Gal(F/M_0)$ such that the twisted $L$-function
$L(E/M_0,\chi,s)$ has a simple zero at $s=1$.
Their theorem requires $N_E, \Delta_{M_0}$ and $p$ to be coprime, but as
they explain this is only necessary to invoke the Gross--Zagier--Zhang formula;
this formula has now proved in complete generality by
Yuan--Zhang--Zhang \cite{YZZ}.

By Tian--Zhang \cite{TZ}, which generalises the earlier work by
Bertolini--Darmon \cite{BD},
the $\chi$-component of $\X$ has multiplicity~1. So $\X$ contains exactly
one copy of the unique $(p\!-\!1)p^n$-dimensional $\Q_p$-irreducible
$p$-adic representation of $\Gal(F/\Q)$.
Its restriction to $H$ is $\rho^{\oplus p^n}$,
and no other representation contributes to $\rho$, so $m_\rho=p^n$ is odd.

(As an alternative to the yet unavailable \cite{YZZ,TZ},
one may bypass the $L$-functions completely by combining
Cornut--Vatsal's \cite{CV} Thm. 4.2 with Nekov\'a\v r's \cite{NekE} Thm. 3.2.
This directly yields a $\chi$ such that the $\chi$-component
of $\X$ has multiplicity~1.)
\end{proof}

\begin{corollary}
For every $E/\Q$, either the Birch--Swinnerton-Dyer rank
formula holds modulo 2 (Conjecture \ref{parityconj}), or
$\sha(E/\Q)$ contains a copy of~$\Q/\Z$.
\end{corollary}


\begin{thebibliography}{29}


\bibitem{BD}
M. Bertolini, H. Darmon,
Iwasawa's Main Conjecture for elliptic curves over
anticyclotomic $\Z_p$-extensions, Annals of Math. 162,
Number 1 (2005), 1--64.

\bibitem{Bir}
B. J. Birch, Conjectures concerning elliptic curves,
Proc. Sympos. Pure Math., Vol. VIII (1965), Amer. Math. Soc., Providence,
R.I, 106--112.

\bibitem{BS}
B. J. Birch, N. M. Stephens,
The parity of the rank of the Mordell-Weil group,
Topology 5 (1966), 295--299.

\bibitem{BFH}
D. Bump, S. Friedberg, and J. Hoffstein, Nonvanishing theorems for
L-functions of modular forms and their derivatives, Invent. Math. 102
(1990), 543--618.

\bibitem{CasIV}
J. W. S. Cassels, Arithmetic on curves of genus 1, IV,
Proof of the Hauptvermutung, J. Reine Angew. Math. 211, (1962), 95--112.

\bibitem{CasVIII}
J. W. S. Cassels,
Arithmetic on curves of genus 1. VIII: On conjectures of Birch and Swinnerton-Dyer,
J. Reine Angew. Math. 217 (1965), 180--199 (1965).

\bibitem{CFKS}
J. Coates, T. Fukaya, K. Kato, R. Sujatha,
Root numbers, Selmer groups and non-commutative Iwasawa theory,
preprint.

\bibitem{CS}
J. Coates, R. Sujatha, appendix to
T. Dokchitser, V. Dokchitser,
Computations in non-commutative Iwasawa theory,
Proc. London Math. Soc. (3) 94 (2006) 211--272.

\bibitem{CV}
C. Cornut, V. Vatsal, Nontriviality of Rankin-Selberg L-functions
and CM points,
in: L-functions and Galois representations (Durham, July 2004),
LMS Lecture Note Series 320 (2007), Cambridge Univ. Press, 121--186.

\bibitem{DN}
C. Diem, N. Naumann, On the structure of Weil restrictions of
Abelian varieties, J.~Ramanujan Math. Soc. 18, No.2 (2003), 153--174.

\bibitem{TV-P}
T. Dokchitser, V. Dokchitser,
Parity of ranks for elliptic curves with a cyclic isogeny,
J. Number Theory 128 (2008), 662--679.

\bibitem{Selfduality}
T. Dokchitser, V. Dokchitser,
Self-duality of Selmer groups, 2007, arxiv: 0705.1899,
to appear in Math. Proc. Cam. Phil. Soc.

\bibitem{Tamroot}
T. Dokchitser, V. Dokchitser,
Regulator constants and the parity conjecture,
2007, arxiv: 0709.2852.

\bibitem{Fis}
T. Fisher, Appendix to V. Dokchitser,
Root numbers of non-abelian twists of elliptic curves,
Proc. London Math. Soc. (3) 91 (2005), 300--324.

\bibitem{Gre}
R. Greenberg, On the Birch and Swinnerton-Dyer conjecture,
Invent. Math. 72, no. 2 (1983), 241--265.

\bibitem{Guo}
L. Guo,
General Selmer groups and critical values of Hecke L-functions,
Math. Ann. 297 no. 2 (1993), 221--233.

\bibitem{HV}
Y. Hachimori, O. Venjakob, Completely faithful Selmer groups over Kummer
extensions, Documenta Mathematica, Extra Volume: Kazuya Kato's Fiftieth
Birthday (2003), 443--478.

\bibitem{Kim}
B. D. Kim,
The Parity Theorem of Elliptic Curves at Primes with Supersingular Reduction,
Compositio Math. 143 (2007) 47--72.

\bibitem{Kob}
S. Kobayashi, The local root number of elliptic curves with wild
ramification, Math. Ann. 323 (2002), 609--623.

\bibitem{Kol}
V. A. Kolyvagin, Euler systems, The Grothendieck Festschrift,
Prog. in Math., Boston, Birkhauser (1990).

\bibitem{Kra}
K. Kramer, Arithmetic of elliptic curves upon quadratic extension,
Trans. Amer. Math. Soc. 264 (1981), 121--135.

\bibitem{KT}
K. Kramer, J. Tunnell, Elliptic curves and local $\epsilon$-factors,
Compositio Math. 46 (1982), 307--352.

\bibitem{MR}
B. Mazur, K. Rubin, Finding large Selmer ranks via an arithmetic theory
of local constants,
Annals of Math. 166 (2), 2007, 579--612.

\bibitem{MilO}
J. S. Milne, On the arithmetic of abelian varieties,
Invent. Math. 17 (1972), 177--190.

\bibitem{MilA}
J. S. Milne, Arithmetic duality theorems,
Perspectives in Mathematics, No. 1, Academic Press, 1986.

\bibitem{Mon}
P. Monsky, Generalizing the Birch--Stephens theorem. I: Modular curves,
Math. Z., 221 (1996), 415--420.

\bibitem{MM}
M. R. Murty and V. K. Murty, Mean values of derivatives of modular
L-series, Annals of Math. 133 (1991), 447--475.

\bibitem{NekS}
J. Nekov\'a\v r, Selmer complexes, Ast\'erisque 310 (2006).

\bibitem{NekE}
J. Nekov\'a\v r,
The Euler system method for CM points on Shimura curves,
in: \hbox{L-functions} and Galois representations (Durham, July 2004),
LMS Lecture Note \hbox{Series} 320 (2007), Cambridge Univ. Press, 471--547.

\bibitem{PS}
B. Poonen, M. Stoll, The Cassels--Tate pairing on polarized abelian varieties,
Annals of Math. 150 (1999), 1109--1149.

\bibitem{RohV}
D. E. Rohrlich, Variation of the root number in families of elliptic curves,
Compositio Math. 87 (1993), 119--151.

\bibitem{SerLi}
J.-P. Serre, Linear Representations of Finite Groups,
GTM 42, Springer Verlag 1977.

\bibitem{Sil2}
J. H. Silverman, Advanced Topics in the Arithmetic of Elliptic Curves,
GTM 151, Springer-Verlag 1994.

\bibitem{TatD}
J. Tate, Duality theorems in Galois cohomology over number fields,
Proc. ICM Stockholm 1962, 234--241.

\bibitem{TatC}
J. Tate, On the conjectures of Birch and Swinnerton-Dyer and a
geometric analog, S\'eminaire Bourbaki, 18e ann\'ee, 1965/66, no. 306.

\bibitem{TZ}
Y. Tian, S. Zhang, Euler system of CM-points on Shimura curves,
in preparation.

\bibitem{Wal}
J.-L. Waldspurger, Correspondences de Shimura et quaternions,
Forum Math. 3 (1991), 219--307.

\bibitem{Yu}
H. Yu, Idempotent relations and the conjecture
of Birch and Swinnerton-Dyer, Math. Ann. 327 (2003), 67--78.

\bibitem{YZZ}
H. Yuan, S. Zhang, W. Zhang, Heights of CM points I: Gross-Zagier formula,
2008, preprint.

\end{thebibliography}
\end{document}